\documentstyle[amscd,amssymb,verbatim,diagrams,12pt]{amsart}
\newcommand{\bH}{{\bf H}}
\newcommand{\KER}{\operatorname{KER}}
\newcommand{\Id}{\operatorname{Id}}

\renewcommand{\mod}{\operatorname{mod}}

\newcommand{\und}{\underline}

\newcommand{\OO}{{\cal O}}

\newcommand{\Lag}{\operatorname{Lag}}

\newcommand{\BB}{{\cal B}}

\newcommand{\SL}{\operatorname{SL}}
\newcommand{\G}{{\Bbb G}}

\newcommand{\Coh}{\operatorname{Coh}}
\newcommand{\GG}{{\cal G}}
\newcommand{\CC}{{\cal C}}

\newcommand{\Spec}{\operatorname{Spec}}

\newcommand{\si}{\sigma}

\newcommand{\Pic}{\operatorname{Pic}}

\newcommand{\ga}{\gamma}
\newcommand{\de}{\delta}

\renewcommand{\ker}{\operatorname{ker}}

\numberwithin{equation}{subsection}

\newtheorem{thm}{Theorem}[subsection]
\newtheorem{prop}[thm]{Proposition}
\newtheorem{lem}[thm]{Lemma}
\newtheorem{cor}[thm]{Corollary}
{  \theoremstyle{definition}
           \newtheorem{defi}[thm]{Definition}
           \newtheorem{rem}[thm]{Remark}
           \newtheorem{rems}[thm]{Remarks}
          \newtheorem{ex}[thm]{Example}
        \newtheorem{exs}[thm]{Examples}
}

\newcommand{\Pf}{\noindent {\it Proof}}
\newcommand{\id}{\operatorname{id}}

\newcommand{\ov}{\overline}

\newcommand{\rk}{\operatorname{rk}}

\newcommand{\Ab}{{\cal A}b}
\newcommand{\FF}{{\cal F}}
\newcommand{\EE}{{\cal E}}
\newcommand{\MM}{{\cal M}}
\newcommand{\TT}{{\cal T}}

\newcommand{\PP}{{\cal P}}

\newcommand{\VV}{{\cal V}}
\newcommand{\SS}{{\cal S}}
\newcommand{\LL}{{\cal L}}
\newcommand{\Du}{{\Bbb D}}

\newcommand{\Hom}{\operatorname{Hom}}

\newcommand{\Ext}{\operatorname{Ext}}
\newcommand{\End}{\operatorname{End}}

\newcommand{\PPic}{{\cal P}ic}

\renewcommand{\a}{\alpha}
\renewcommand{\b}{\beta}

\newcommand{\De}{\Delta}
\newcommand{\bj}{{\bf j}}
\newcommand{\la}{\lambda}

\newcommand{\N}{{\Bbb N}}
\newcommand{\R}{{\Bbb R}}
\newcommand{\Z}{{\Bbb Z}}
\newcommand{\Q}{{\Bbb Q}}
\newcommand{\bbT}{{\Bbb T}}
\newcommand{\La}{\Lambda}
\newcommand{\Ga}{\Gamma}

\newcommand{\wt}{\widetilde}
\newcommand{\ot}{\otimes}

\newcommand{\sub}{\subset}
\newcommand{\ed}{\qed\vspace{3mm}}

\newcommand{\rep}{\operatorname{rep}}
\newcommand{\fib}{\operatorname{fib}}
\newcommand{\Fun}{\operatorname{Fun}}
\newcommand{\Nm}{\operatorname{Nm}}
\newcommand{\cha}{\operatorname{char}}

\newcommand{\Br}{\operatorname{Br}}

\title{Lagrangian-invariant sheaves and functors for abelian varieties}
\author{Alexander Polishchuk}
\address{Department of Mathematics, University of Oregon, Eugene, OR 97405}
\email{apolish@@uoregon.edu}
\thanks{Supported in part by the NSF grant DMS-1001364}

\begin{document}
\begin{abstract} We partially generalize the theory of semihomogeneous bundles on an abelian
variety $A$ developed by Mukai~\cite{Mukai}. This involves considering abelian subvarieties
$Y\sub X_A=A\times\hat{A}$ and studying coherent sheaves on $A$ invariant under
the action of $Y$. The natural condition to impose on $Y$ is that of
being {\it Lagrangian} with respect to a certain skew-symmetric biextension $\EE$ of $X_A\times X_A$
by $\G_m$. We prove that in this case any $Y$-invariant sheaf is a direct sum of several copies 
of a single coherent sheaf. We call such sheaves {\it Lagrangian-invariant} (or {\it LI-sheaves}). 
We also study {\it LI-functors} $D^b(A)\to D^b(B)$ 
associated with kernels in $D^b(A\times B)$ that are invariant with respect to some Lagrangian
subvariety in $X_A\times X_B$. We calculate their composition and prove that in characteristic zero it can be decomposed into a direct sum of LI-functors. In the case $B=A$ this
leads to an interesting central extension of the group of symplectic automorphisms of $X_A$ in the category of abelian varieties up to isogeny.
\end{abstract}
\maketitle

\bigskip

\centerline{\sc Introduction}

Derived categories of coherent sheaves on smooth projective varieties have been playing an increasingly significant role in modern algebraic geometry. They serve as a motivation and
a testing ground for constructions in noncommutative geometry of smooth proper dg-categories (see \cite{KS}) and have interesting links to representation theory, mirror symmetry and other areas (see 
\cite{Bridge-ICM} and references therein).
One of the first results that revealed the importance of derived categories
was the equivalence of the derived categories 
$D^b(A)\simeq D^b(\hat{A})$ for an abelian variety $A$ and its dual $\hat{A}$.
established by Mukai~\cite{Muk-FM}. 
The present paper continues the line of study of \cite{O}, where it was shown
that the group of autoequivalences of $D^b(A)$ is closely related to the group $U(A\times\hat{A})$ of
``symplectic" automorphisms of $A\times\hat{A}$ (the term ``symplectic" is explained below). 
Namely, we study 
endofunctors of $D^b(A)$ related to symplectic automorphisms of $A\times\hat{A}$ 
in the bigger category $\Ab_{\Q}$ of abelian varieties up to isogeny.
More generally, we study a natural class of {\it LI-functors} between
derived categories of abelian varieties 
governed by {\it Lagrangian correspondences} between the corresponding
varieties $A\times\hat{A}$. 
Our motivation is twofold. On the one hand, such functors could be useful
for the Manin's real multiplication program (see \cite{Manin}) if one tries to approach it using
the categories of holomorphic bundles on noncommutative tori described in \cite{PS}.
On the other hand, as will be shown in a sequel to this paper, using these functors one can
realize the Bridgeland stability space of an abelian surface $A$ as a homogeneous space for certain 
covering of the group $U(A\times\hat{A},\R)$ of $\R$-points of the algebraic group of
symplectic automorphisms of $A\times\hat{A}$ (more precisely, one should consider the connected component of the stability space described in \cite{Bridge2}).

As in \cite{P-sym}, with every abelian variety $A$ we associate a ``symplectic" object: the abelian
variety $X_A=A\times\hat{A}$ together with a skew-symmetric biextension
$\EE=p_{14}^*\PP\ot p_{23}^*\PP^{-1}$ of $X_A\times X_A$ (where $\PP$ is the Poincar\'e line bundle
on $A\times\hat{A}$). There is an analog of the Schr\"odinger
representation of the Heisenberg group in this setting obtained by considering the action of a certain
Heisenberg extension $\bH$ of
$X_A$ on $D^b(A)$ by translations and tensoring with line bundles in $\Pic^0(A)$ 
(see Sec.\ \ref{sym-sec}).
One can show that under some technical assumptions,
the stabilizer of an endosimple object of $D^b(A)$ is always an isotropic
subgroup $Y\sub X_A$ with respect to the biextension $\EE$ (see Remark \ref{moduli-rem}). 
We study the case when this subgroup is a {\it Lagrangian}, i.e., when we have a
duality between $Y$ and $X_A/Y$ (we show that this is equivalent to $Y$ having the maximal
possible dimension $\dim A=\dim X_A/2$). 
Our first main result is the description of the subcategory of invariants in $D^b(A)$
of a Lagrangian abelian subvariety $Y\sub X_A$ (see Theorem \ref{lag-inv-thm}). 
More precisely, to define this subcategory of invariants we need to choose a line bundle $\a$ on $Y$ 
that gives a lifting of $Y$ to the Heisenberg groupoid $\bH$ of $A$ (see Sec.\ \ref{isotr-sec}). Then
we can define the subcategory of $(Y,\a)$-invariant objects in $D^b(A)$. We show that 
this subcategory is generated by a single endosimple coherent sheaf $S_{Y,\a}$ and that every
$(Y,\a)$-invariant coherent sheaf is a direct sum of several copies of $S_{Y,\a}$. We call such sheaves
{\it Lagrangian-invariant} ({\it LI-sheaves} for short).
Note that the case when the projection $Y\to A$ is an isogeny corresponds to
semihomogeneous vector bundles on $A$ (see \cite{Mukai}). However, our condition of 
$(Y,\a)$-invariance is stronger than the condition of invariance under all
functors corresponding to points in $Y(k)$. Namely, the condition for a sheaf $F\in\Coh(A)$
to be $(Y,\a)$-invariant is given by a certain isomorphism over $A\times Y$ (see Sec.\ \ref{inv-sec}).
This accounts for the fact that nontrivial extensions of $S_{Y,\a}$ with itself
are not $(Y,\a)$-invariant (see Example \ref{homog-ex}).

Note that finite Heisenberg groups appear in our picture as well. Namely,
for a pair of Lagrangian subvarieties $(Y,Z)$ in $X_A$ equipped with liftings to the Heisenberg groupoid, we have a natural group extension $G$ of the intersection $Y\cap Z$  by $\G_m$, which is a Heisenberg group in the case when $Y\cap Z$ is finite (this construction goes back to
\cite{P-sym}). 
Generalizing Mumford's theory of theta groups we show that
the space $\Hom^*_{D^b(A)}(S_{Y,\a}, S_{Z,\b})$ is an
irreducible representation of the Heisenberg group $G$ (see Corollary \ref{Hom-cor}).

For a pair of abelian varieties $A$ and $B$ we consider Lagrangian correspondences 
from $X_A$ to $X_B$, i.e., Lagrangian abelian subvarieties $L\sub X_A\times X_B$ with respect
to the difference of biextensions $\EE$ pulled-back from $X_A\times X_A$ and $X_B\times X_B$.
More precisely, we consider pairs $(L,\a)$, where $\a$ is a line bundle on $L$ that gives a lifting of $L$
to the Heisenberg groupoid of $X_A\times X_B$. As above, with each pair $(L,\a)$ we associate
the LI-sheaf $S_{L,\a}$ on $A\times B$. Then we define the corresponding {\it LI-functor} 
$\Phi_{L,\a}:D^b(A)\to D^b(B)$ as the Fourier-Mukai functor associated with $S_{L,\a}$. 
Our main result is the calculation of the composition of functors of this form.
We do this in two steps. First, we introduce {\it generalized Lagrangians} in $X_A$
by allowing for $Y$ not to be connected and for the homomorphism $Y\to X_A$ to have finite kernel.
The condition for $Y$ to be Lagrangian is formulated in terms of duality between $Y$ and
the complex $[Y\to X_A]$. The context for such a duality is provided by the category of
{\it orbi-abelian varieties} studied in \cite[Sec.\ 4]{P-FM}. The theory of invariants of Lagrangians in
$D^b(A)$ carries over to generalized Lagrangians. We define composition of generalized Lagrangian
correspondences (which is again generalized Lagrangian) and
show that the composition of LI-functors is compatible with the composition 
of correspondences up to a shift (see Theorem \ref{Lag-conv-thm}). 
The second step consists of the study of the relation between the LI-sheaf $S_{Y,\a}$ associated with 
a generalized Lagrangian $(Y,\a)$ for $X_A$ and the corresponding usual Lagrangian obtained by passing to the image in $X_A$
of the connected component of zero $Y_0\sub Y$ (see Proposition \ref{lag-inv-prop}).
The final result (in characteristic zero) is that the composition of two LI-functors decomposes as an explicit direct sum of LI-functors. In the case $B=A$ we obtain an explicit description of
a certain semiring of endofunctors of $D^b(A)$ in terms of the group $U(X_A,\Q)$ of symplectic automorphisms of $X_A$ in $\Ab_\Q$.
Considering these endofunctors up to a natural action of $X_A(k)$ we get a simpler
picture: the resulting monoid has a group completion that can be identified with an explicit
central extension of $U(X_A,\Q)$ by $\Z\times\Q^*$. The $\Z$-part of this central extension is described by essentially the same
$2$-cocycle as in the picture with autoequivalences of $D^b(A)$ (see \cite[Sec.\ 4]{O}). The $\Q^*$-part
of the central extension is a new feature of the picture with endofunctors,
and we show that it can be nontrivial in the case of elliptic curves
with complex multiplication.

The paper is organized as follows. In Section \ref{prelim-sec} we gather necessary results
about finite Heisenberg group schemes and twisted equivariant sheaves, along with
some other background facts. In Section \ref{Heis-sec} we review the symplectic setting from 
\cite{P-sym}. The main object here is the Heisenberg groupoid $\bH$ associated with a symplectically
self-dual abelian variety (equipped with some extra data). We define generalized isotropic
(resp., Lagrangian) pairs and the corresponding representations of $\bH$. The main result of this section
is Theorem \ref{lag-inv-thm} describing the subcategories of invariants of a generalized Lagrangian pair  
in standard representations of $\bH$. Another important technical result is Proposition
\ref{lag-inv-prop} giving a relation between invariants for a generalized Lagrangian $i:Y\to X$
and the associated Lagrangian $i(Y_0)\sub X$. In Section \ref{corr-sec} we study LI-kernels
and LI-functors. First, we define generalized Lagrangian correspondences 
(or {\it g-Lag-correspondences}) and show
that they can be composed under some assumptions. Then in \ref{ker-fun-sec}
we prove the compatibility
between the composition of g-Lag-correspondences and the composition of the corresponding
LI-functors up to a shift (see Theorems \ref{main-conv-thm} and Theorem \ref{Lag-conv-thm}).
Finally in \ref{central-ext-sec} assuming that the characteristic is zero we describe
a semiring of endofunctors of $D^b(A)$ built from LI-functors and the related central
extension of the group $U(X_A,\Q)$.


\medskip

\noindent
{\it Notation and conventions}: We consider schemes over a fixed algebraically closed ground field $k$. 
We say that an object $F$ of a $k$-linear category $\CC$ is {\it endosimple} if
$\Hom_\CC(F,F)=k$.
For a scheme $X$ we denote by $D^b(X)$ the bounded derived category of coherent sheaves on $X$.
We say that an object $F\in D^b(X)$ is {\it cohomologically pure} if there exists a coherent sheaf $H$
such that $F\simeq H[n]$ for $n\in\Z$. For $F_1,F_2\in D^b(X)$ we set
$\Hom^*_{D^b(X)}(F_1,F_2)=\bigoplus_{i}\Hom_{D^b(X)}(F_1,F_2[i])$.
We denote by $\Ab_\Q$ the category of abelian varieties up to an isogeny (i.e., the localization
of the category of abelian varieties over $k$ with respect to the class of isogenies).
For a finite group scheme $G$ by a representation of $G$ on a vector space 
we always mean a finite-dimensional representation.
For a group scheme $G$ we denote by $G_0$ the connected component of the neutral element
in $G$ with the reduced scheme structure, and we set $\pi_0(G)=G/G_0$.
All our biextensions are biextensions by $\G_m$. Most of the facts about biextensions we use
are summarized in \cite[II.10]{P-ab}. 
We always normalize the Poincar\'e line bundle $\PP=\PP_A$ on $A\times\hat{A}$ by requiring it
to be trivial over $\{0\}\times\hat{A}$ and $A\times\{0\}$ (thus making it a biextension of
$A\times\hat{A}$). For a line bundle $\LL$ over a commutative group scheme $X$ we denote
by $\La(\LL)$ the line bundle on $X\times X$ given by $\La(\LL)_{x,x'}=\LL_{x+x'}\ot\LL^{-1}_x\ot
\LL^{-1}_{x'}$.

\section{Preliminaries}\label{prelim-sec}

\subsection{Finite Heisenberg group schemes}\label{Heis-gr-sec}

In this section we recall some facts about Heisenberg group schemes and their representations
of weight one following \cite[Sec.\ 12.2]{P-ab} (see also \cite[Ch.\ V]{MB} and \cite[Sec.\ 1,2]{P-maslov}).

A {\it finite Heisenberg group scheme} 
$G$ is a central extension of a finite commutative group scheme $K$ by $\G_m$,
such that the commutator pairing $K\times K\to \G_m$ is perfect, i.e., induces an isomorphism
$K\simeq K^*$, where $K^*$ is the Cartier dual to $K$.

In general, when $G$ is a central extension of some group scheme $K$ by $\G_m$
we say that a finite-dimensional representation $V$ of $G$ is of {\it weight one} if
$\G_m\sub G$ acts on $V$ via the identity character.
All weight-$1$ representations of a finite Heisenberg group scheme
are direct sums of several copies of the unique
irreducible representation of weight one, called {\it Schr\"odinger representation}.
Its dimension is equal to $|K|^{1/2}$, where $|K|$ is the order of $K$.

Sometimes it is convenient to describe a Heisenberg extension $G$ in terms of the line bundle $L$
over $K$ associated with $G$ viewed as a $\G_m$-torsor over $K$. The group law on $G$
corresponds to an isomorphism
\begin{equation}\label{Heis-tors-isom}
L_{x+y}\simeq L_x\ot L_y.
\end{equation}
of line bundles on $K\times K$. 
A representation of $G$ of weight one on a finite dimensional vector space $V$ can
be described as a morphism
$L\ot V\to \OO_K\ot V$ of bundles over $K$ satisfying the natural
compatibility with \eqref{Heis-tors-isom}. We denote by $G-\rep_1$ the category of
weight-$1$ representations of $G$.

Any subgroup $I\sub K$, isotropic with respect to the commutator pairing, can be lifted to a subgroup
of $G$. Equivalently, there is a trivialization of $L$ over $I$, compatible with
\eqref{Heis-tors-isom}. Any maximal isotropic subgroup of $K$ has order $|K|^{1/2}$.
The Schr\"odinger representation can be realized explicitly as follows.
The space $H^0(K,L)$ has natural commuting left and right $G$-actions of weight one:
the corresponding isomorphisms 
$$L\ot H^0(K,L)\simeq \OO_K\ot H^0(K,L)$$ 
are obtained
from \eqref{Heis-tors-isom} by taking the push-forward with respect to the projection $K\times K\to K$
to one of the factors. Now if $I\sub K$ is a maximal isotropic subgroup then we can lift it
to a subgroup of $G$ and take invariants $H^0(K,L)^I$ with respect to the right action.
This becomes an irreducible $G$-representation of weight one. Trivializing the line
bundle $L$ over $K$ we can identify $H^0(K,L)^I$ with the subspace of functions $f$ on $K$ 
satisfying certain quasi-periodicity with respect to translations by $I$,
which is the customary description of the Schr\"odinger representation.
If $V$ is a Schr\"odinger representation of $G$ and $I\sub K$ is a maximal isotropic subgroup
with a lifting to $G$ then $V^I$ is one-dimensional. Hence, 
if $W$ is any weight-$1$ representation of $G$ then we have an isomorphism of $G$-representations
(depending on a trivialization of $V^I$):
$$W\simeq V\ot W^I.$$

Note that if $V$ is a Schr\"odinger representation of $G$ then $V^*$ is a Schr\"odinger representation
of $G^{op}$ (using the isomorphism $G\to G^{op}:g\mapsto g^{-1}$ we can also view $V^*$ as
a representation of $G$ of weight $-1$).
For the entire space $H^0(K,L)$ we have an isomorphism of $G\times G^{op}$-modules
\begin{equation}\label{End-eq}
H^0(K,L)\simeq V\ot V^*. 
\end{equation}

Often when we talk about an isotropic subgroup in $K$ lifted to $G$ we just say that
$I$ is an isotropic subgroup in $G$. If $I$ is such a subgroup then the normalizer subgroup $N(I)$
contains $\G_m$ and $N(I)/\G_m=I^{\perp}$, the orthogonal complement to $I$ with respect
to the commutator pairing. The quotient $N(I)/I$ is a Heisenberg extension of $I^{\perp}/I$.
If $V$ is a Schr\"odinger representation of $G$ 
then the space $V^I$ of $I$-invariants is a Schr\"odinger representation of $N(I)/I$.
Viewed as a representation of $I$, 
$V$ is a direct sum of $\dim(V^I)$ copies of the regular representation 
$H^0(I,\OO)^*\simeq H^0(I^*,\OO)$.
 Hence, in the case when $I^*$ is reduced we get an isomorphism 
$$V\simeq\bigoplus_{\chi\in I^*}V^I\ot \chi$$
of $I$-representations.

\subsection{1-cocycles with values in Picard stacks and twisted equivariant sheaves}\label{twisted-sec}

Let us recall briefly some definitions from \cite[Sec.\ 3.4]{P-FM}.
Let $X$ be a scheme, $G$ a group scheme acting on $X$. 
A {\it $1$-cocycle of $G$ with values in the Picard
stack $\PPic(X)$ of line bundles on $X$} is a line bundle $L$ on $G\times X$ equipped
with an isomorphism
\begin{equation}\label{1-coc-def-eq}
L_{g_1g_2,x}\simeq L_{g_1,g_2x}\ot L_{g_2,x}
\end{equation}
on $G\times G\times X$
such that the following diagram of isomorphisms on $G\times G\times G\times X$ is commutatve:
\begin{equation}\label{1-coc-comm-eq}
\begin{diagram}
L_{g_1g_2g_3,x}&\rTo{}&
L_{g_1,g_2g_3x}\ot L_{g_2g_3,x}\\
\dTo{}&&\dTo{}\\
L_{g_1g_2,g_3x}\ot L_{g_3,x}&\rTo{}&L_{g_1,g_2g_3x}\ot L_{g_2,g_3x}\ot L_{g_3,x}
\end{diagram}
\end{equation}
(cf. \cite[(3.4.2)]{P-FM}). Given such a $1$-cocycle $L$
we can consider $L$-twisted $G$-equivariant sheaves on $X$ as sheaves equipped
with an isomorphism
\begin{equation}\label{L-twisted-eq}
F_{gx}\simeq L_{g,x}\ot F_x
\end{equation}
on $G\times X$ subject to the commutativity of the diagram
\begin{equation}\label{L-twisted-comm-eq}
\begin{diagram}
F_{g_1g_2x}&\rTo{}&
L_{g_1,g_2x}\ot F_{g_2x}\\
\dTo{}&&\dTo{}\\
L_{g_1g_2,x}\ot F_{x}&\rTo{}&L_{g_1,g_2x}\ot L_{g_2,x}\ot F_{x}
\end{diagram}
\end{equation}
on $G\times G\times X$ (cf. \cite[Sec.\ 3.4]{P-FM}). We denote the category of
$L$-twisted $G$-equivariant coherent sheaves on $X$ by $\Coh^L_H(X)$.
Similarly, we define the category $D^L_H(X)$ of
$L$-twisted $G$-equivariant objects in the derived category $D^b(X)$
as objects $F\in D^b(X)$ equipped with an isomorphism \eqref{L-twisted-eq} in
$D^b(G\times X)$, subject to the same compatibilities.

For example, if $G\to K$ is Heisenberg group extension then the underlying $\G_m$-torsor
$L$ over $K$ has a structure of a $1$-cocycle of $K$ with values
in the Picard stack $\PPic(k)$ of $1$-dimensional vector spaces, and the category $G-\rep_1$
can be identified with the category of $L$-twisted $K$-equivariant coherent sheaves on
$\Spec(k)$.

\begin{lem}\label{center-lem}
Let $G$ be a central extension of a finite commutative group scheme
$K$ by $\G_m$ with the center $C\times\G_m$, such that
the induced central extension $\ov{G}=G/C\to K/C$ is a Heisenberg group. 
Let $p_{\ov{G}}:G\to \ov{G}$ and $p_K:G\to K$ be the projections, and let
$V$ be a Schr\"odinger representation of $\ov{G}$.
Then there exists an equivalence
\begin{equation}\label{Phi-eq}
\Phi: G-\rep_1\to C-\rep
\end{equation}
such that $\Phi(p_{\ov{G}}^*V)$ is the trivial one-dimensional representation of $C$, and
for a character $\chi:K\to\G_m$ we have
\begin{equation}\label{kappa-comp-eq}
\Phi(\kappa_{\chi}^*W)\simeq \chi|_C\ot \Phi(W),
\end{equation}
where $\kappa_{\chi}:G\to G$ is the automorphism sending $g$ to 
$\chi(p_K(g))\cdot g$.
\end{lem}

\Pf . By \cite[Lem.\ 3.5.3(i)]{P-FM}, the extension of commutative groups
\begin{equation}\label{C-K-ext-eq}
0\to C\to K\rTo{\pi} K/C\to 0
\end{equation}
gives rise to a biextension $\BB$ of $K/C\times C^*$. Furthermore, if
$\wt{\MM}$ is the $1$-cocycle of $K$ with values in $\PPic(k)$ corresponding to the extension
$G\to K$, then we have $\wt{\MM}\simeq\pi^*\MM_0$, where $\MM_0$ is the $1$-cocycle
of $K/C$ with values in $\PPic(k)$ corresponding to the extension $\ov{G}\to K/C$.
We can extend $\MM_0$ to a $1$-cocycle $\MM$ of $K/C$ 
with values in $\PPic(C^*)$ by the formula
$$\MM_{x,c^*}=\BB_{x,c^*}\ot (\MM_0)_x,$$
where $x\in K/C$, $c^*\in C^*$. 
Then one has a natural isomorphism of $k$-algebras
$$\alpha: p_*\OO_{C^*}[K/C]^{\MM}\rTo{\sim} k[K]^{\wt{\MM}},$$
where $p:C^*\to\Spec(k)$ is the projection to the point 
(see \cite[Prop.\ 3.5.4]{P-FM} for a more general result). 
Here for a group $H$ acting trivially on a scheme $X$ and a $1$-cocycle $\LL$ of $H$ with
values in $\PPic(X)$ we denote by $\OO_X[H]^{\LL}$ the twisted group algebra of
$H$ over $X$ (see \cite[(3.4.9)]{P-FM}).
The equivalence between the categories of modules induced by $\a$ can be interpreted as an 
equivalence
$$\Coh^{\MM}_{K/C}(C^*)\simeq\Coh^{\wt{\MM}}_K(\Spec(k)).$$
As we have seen above, the category on the right can be identified with 
$G-\rep_1$. On the other hand, $\Coh^{\MM}_{K/C}(C^*)$ is equivalent to the category of
weight-$1$ representations of the finite Heisenberg group scheme $\ov{\GG}$ over $C^*$,
where $\ov{\GG}$ is the extension of the constant finite group scheme $(K/C)_{C^*}$
by $\G_m$ associated with the $1$-cocycle $\MM$. Since the obstacle to the existence of
a Schr\"odinger representation is an element of the Brauer group $\Br(C^*)=0$
(see \cite[Sec.\ 2]{P-maslov}), the latter category is equivalent to 
$\Coh(C^*)\simeq C-\rep$.
Thus, we obtain the required equivalence \eqref{Phi-eq}.

It is easy to check that $\alpha$ is compatible with natural homomorphisms of both sides
to $k[K/C]^{\MM_0}$. Note that modules over $k[K/C]^{\MM_0}$ are exactly weight-$1$ representations
of $\ov{G}$. This compatibility implies that $p_{\ov{G}}^*V$ corresponds to the trivial
representation of $C$ under the equivalence $\Phi$.
The compatibility \eqref{kappa-comp-eq} 
with automorphisms $\kappa_{\chi}:G\to G$ follows from the
fact that the whole construction is functorial with respect to isomorphisms $G\to G'$, identical
on $\G_m$.
\ed

We will use the following description of the categories of 
twisted equivariant sheaves associated with some transitive actions.

\begin{prop}\label{cocycle-prop} 
Let $H$ be a group scheme acting transitively on $X$ with finite stabilizers, and let $L_{h,x}$ be a $1$-cocycle for the action of $H$ on $X$. 

\noindent (i)
For a point $x\in X$ there is a natural central extension 
$$1\to \G_m\to G_x\to H_x\to 1$$ 
of the stabilizer subgroup $H_x$ of
$x$ with the underlying line bundle $L|_{H_x\times\{x\}}$.
Consider the line bundle $L(x):=L|_{H\times\{x\}}$ over $H$, and
let $\pi(x):H\to X$ be the map sending $h$ to $hx$.
Then $L(x)$ has a natural
$(\id_H\times\pi(x))^*L$-twisted $H$-equivariant structure with respect to the
regular left action of $H$ on itself.
The right action of $H_x$ on $H$ extends to an action of $G_x$ on $L(x)$
respecting the above twisted $H$-equivariant structure.

The map $F\mapsto F|_x$ extends to an equivalence
\begin{equation}\label{cocycle-eq}
\fib_x: \Coh^L_H(X)\to G_x-\rep_1
\end{equation}
of $L$-twisted $H$-equivariant coherent sheaves on $X$ with the category of weight-$1$ representations of $G_x$. 

\noindent (ii)
For an element $h\in H$ the isomorphism $H_x\to H_{hx}$ given 
by the conjugation by $h$ extends to an isomorphism of central extensions 
\begin{equation}\label{central-ext-isom}
\kappa_h:G_x\to G_{hx}.
\end{equation}
 Furthermore, there is a natural isomorphism of functors
\begin{equation}\label{conj-fib-isom}
L_{h,x}\ot\fib_x\simeq \kappa_{h}^*\circ \fib_{hx},
\end{equation}
where $\kappa_{h}^*:G_{gx}-\rep_1\to G_x-\rep_1$
is the equivalence induced by $\kappa_h$ ($L_{h,x}$ is just a one-dimensional vector space
one has to insert to make the isomorphism canonical).

\noindent (iii)
Now  let $f:X\to X$ be an automorphism commuting with the $H$-action. Let $M$ be a line bundle on $X$ together with an isomorphism of cocycles
$$\tau_f:M_x\otimes M_{hx}^{-1}\otimes L_{h,x}\to L_{h,f(x)}$$
(so $M$ is a coboundary for $L\ot (\id\times f)^*L^{-1}$).
Then the map $F\mapsto f^*F\ot M$ extends to an  autoequivalence $\iota_{f,M}$ of $\Coh^L_H(X)$, so that the following diagram of functors is commutative up to an isomorphism
\begin{equation}\label{fib-aut-diagram}
\begin{diagram}
\Coh^L_H(X)&\rTo{\fib_x}&G_x-\rep_1\\
\dTo{\iota_{f,M}}&&\dTo{(\alpha_f)^*}\\
\Coh^L_H(X)&\rTo{\fib_x}&G_x-\rep_1
\end{diagram}
\end{equation}
where $\alpha_f$
is the composition
$$G_x\rTo{\ov{\tau}_f} G_{f(x)} \rTo{\kappa_{h_0}^{-1}} G_x$$
for some $h_0\in H$ such that $h_0x=f(x)$.
Here the homomorphism $\ov{\tau}_f$ is given by the restriction of
$\tau_f$ to $H_x\times\{x\}$.
\end{prop}

\Pf . (i) Since the embedding of the point $x$ into $X$ is an
 $H_x$-morphism, the restriction of $L|_{H_x\times \{x\}}$
has a natural structure of a $1$-cocycle of $H_x$ with values in the Picard groupoid of the point, 
hence it defines a central extension of  $H_x$ by $\G_m$. 
The $(\id_H\times\pi(x))^*L$-twisted $H$-equivariant structure on $L(x)$ is
given by the isomorphism
$$L(x)_{h'h}\simeq L_{h',\pi(x)(h)}\ot L(x)_{h},$$
where $h,h'\in H$,
obtained from \eqref{1-coc-def-eq} by the restriction to $H\times H\times\{x\}$.
Restricting \eqref{1-coc-def-eq} further to $H\times H_x\times\{x\}$ we obtain an isomorphism
\begin{equation}\label{L(x)-right-action-eq}
L(x)_{hh_0}\simeq L(x)_h\ot L_{h_0,x},
\end{equation}
where $h\in H$, $h_0\in H_x$. This isomorphism can be interpreted as a right action
of $G_x$ on $L(x)$. The fact that this action commutes with the twisted left $H$-equivariant
structure on $L(x)$ follows from the commutativity of diagram \eqref{1-coc-comm-eq}.

For an $L$-twisted $H$-equivariant coherent sheaf $F$ on $X$ we can restrict the
structure isomorphism
$$F_{hx}\simeq L_{h,x}\ot F_x$$
to $H_x\times\{x\}\sub H\times X$ and get an isomorphism
\begin{equation}\label{a-h0-x-eq}
F_x\rTo{a(h_0,x)} L_{h_0,x}\ot F_x
\end{equation}
for $h_0\in H_x$
satisfying the natural compatibility condition on $H_x\times H_x$. Thus, we can view
$F_x$ as a weight-$1$ representation of $G_x$.
This gives the functor $\fib_x$.

Conversely, starting with a representation $W$ of $G_x$ of weight $1$ we can construct an $L$-twisted $H$-equivariant sheaf on $X$ as follows. Consider
the sheaf $\wt{F}=W\ot L(x)$ on $H$. Viewing $W$ as a right module over $G_x$ of weight $-1$ 
(via the inversion map on $G_x$) and taking the tensor product of this action with the right
action of $G_x$ on $L(x)$, we obtain a right action of $H_x$ on $\wt{F}$.
Hence, $\wt{F}$ descends to a sheaf on $H/H_x\simeq X$. Now the 
$(\id_H\times\pi(x))^*L$-twisted $H$-equivariant structure on $L(x)$
induces an $L$-twisted $H$-equivariant structure on $F$. 

\noindent (ii) The cocycle structure on $L$ gives an isomorphism
$$L_{hh_0h^{-1},hx}\ot L_{h,x}\simeq L_{hh_0,x}$$
for $h, h_0\in H$. Combining it with the isomorphism
$$L_{h,h_0x}\ot L_{h_0,x}\simeq L_{hh_0,x}$$
we get an isomorphism 
$$L_{hh_0h^{-1},hx}\ot (L_{h,x}\ot L_{h,h_0x}^{-1})\simeq L_{h_0,x}$$
for fixed $h\in H$, compatible with the 
 cocycle structures in $(h_0,x)$.
The restriction to $h_0\in H_x$ gives an isomorphism
\begin{equation}\label{L-h-conj-isom}
L_{hh_0h^{-1},hx}\simeq L_{h_0,x},
\end{equation}
compatible with the central extension structures, so we obtain
the required isomorhism \eqref{central-ext-isom}.

Now \eqref{conj-fib-isom} is given by the isomorphism
$$\phi: L_{h,x}\ot F_x\rTo{\sim} F_{hx}$$
which is a part of the $L$-twisted $H$-equivariant structure on $F$. To check
the compatibility of $\phi$ with the $G_x$-action via the isomorphism \eqref{central-ext-isom}, let us
apply the commutativity of \eqref{L-twisted-comm-eq} for the pairs $(h,h_0)$
and $(hh_0h^{-1},h)$, where $h_0\in H_x$. We get the following commutative diagram
(in which we skipped the tensor product signs for brevity)
\begin{diagram}
L_{h,h_0x} F_{h_0x} &\rTo{\phi^{-1}}& F_{hh_0x} &\rTo^=&F_{(hh_0h^{-1})hx}&
\rTo{a(hh_0h^{-1},hx)}&L_{hh_0h^{-1},hx} F_{hx}\\
\dTo{a(h_0,x)}&&\dTo{}&&\dTo{}&&\dTo{\phi}\\
L_{h,h_0x} L_{h_0,x} F_x&\rTo{}& L_{hh_0,x} F_x&\rTo^=&
L_{(hh_0h^{-1})h,x} F_x&\rTo{}& L_{hh_0h^{-1},hx} L_{h,x} F_{x}
\end{diagram}
where $a(\cdot,\cdot)$ are the maps \eqref{a-h0-x-eq} inducing the action of $G_x$ (resp., $G_{hx}$)
on $F_x$ (resp., $F_{hx}$). The required compatibility follows from the fact that
the composition of arrows in the bottom row corresponds to the isomorphism \eqref{L-h-conj-isom}.

\noindent (iii)
The functoriality of the construction gives the commutative diagram
\begin{diagram}
\Coh^L_H(X)&\rTo{\fib_{f(x)}}&G_{f(x)}-\rep_1\\
\dTo{\iota_{f,M}}&&\dTo{\ov{\tau}_f^*}\\
\Coh^L_H(X)&\rTo{\fib_x}&G_x-\rep_1
\end{diagram}
The diagram \eqref{fib-aut-diagram} follows by applying the isomorphism of functors
established in (ii).
\ed

\begin{lem}\label{push-forward-lem} In the situation of Proposition \ref{cocycle-prop}(i)
the bundle 
$$\VV=\pi(x)_*(L|_{H\times\{x\}})$$ on $X$ has a natural 
$L$-twisted $H$-equivariant structure and a commuting right $G_x$-action (where $G_x$ acts trivially
on $X$). Assume in addition that $H_x$ is commutative and $G_x$ is its Heisenberg extension.
Let $I\sub G_x$ be a maximal isotropic subgroup. Then the bundle $\VV^I$ of $I$-invariants is 
a simple object in $\Coh^L_H(X)$ (i.e., $\fib_x(\VV^I)$ is a Schr\"odinger representation).
Also, one has the following isomorphism of sheaves on $H$ compatible with the right $H_x$-action:
\begin{equation}\label{pi(x)-pull-back-eq}
\pi(x)^*\VV^I\simeq V\ot L|_{H\times\{x\}},
\end{equation}
where $V$ is a Schr\"odinger representation of $G_x$, viewed as a representation of
$G_x^{op}$ of weight $-1$.
\end{lem}

\Pf . The $L$-twisted $H$-equivariant structure and the right $G_x$-action on 
$\VV=\pi(x)_*(L(x))$ are
induced by the $(\id_H\times\pi(x))^*L$-twisted $H$-equivariant structure and the right $G_x$-action
on $L(x)$, respectively.
In the case when $G_x$ is a Heisenberg extension of $H_x$ and $I\sub G_x$ is a maximal isotropic
subgroup, we have an isomorphism
$$\fib_x(\VV^I)\simeq H^0(H_x, L|_{H_x\times\{x\}})^I$$
compatible with the left $G_x$-action. But the right-hand side is a standard model for the
Schr\"odinger representation of $G$ (see Sec.\ \ref{Heis-gr-sec}). Hence, $\VV^I$ is simple
in $\Coh^L_H(X)$.

Let $V$ be a Schr\"odinger representation of $G_x$. Viewing $V^*$ as a Schr\"odinger representation
of $G_x^{op}$,  we obtain
an isomorphism
$$\VV\simeq V^*\ot\VV^I,$$
compatible with $L$-twisted $H$-equivariant structures and the right $G_x$-actions.
Note also that we have a natural isomorphism $H\times_X H\simeq H\times H_x$
(where in the first product we use the map $\pi(x):H\to X$). This gives an isomorphism
$$\pi(x)^*\VV=\pi(x)^*\pi(x)_*(L(x))\simeq p_{H*}m^*(L(x)),$$
where $m:H\times H_x\to H$ is given by the group law and $p_H:H\times H_x\to H$
is the projection. Taking into account the isomorphism \eqref{L(x)-right-action-eq}
we get 
$$\pi(x)^*\VV=H^0(H_x, L(x)|_{H_x})\ot L(x).$$
Hence, we obtain an isomorphism
$$V^*\ot\pi(x)^*\VV^I\simeq \pi(x)^*\VV\simeq H^0(H_x,L(x)|_{H_x})\ot L(x)\simeq
V^*\ot V\ot L(x)$$
compatible with $G_x\times G_x^{op}$-actions, which leads to \eqref{pi(x)-pull-back-eq}.
\ed

\begin{ex} Let $A$ be an abelian variety, $\LL$ a line bundle on $A$ (trivialized at zero)
$K\sub A$ a finite subgroup such that $\phi_\LL(K)=0$, where 
$\phi_\LL:A\to\hat{A}$ is the homomorphism associated with $A$.
Consider the biextension $\La(\LL)=m^*\LL\ot p_1^*\LL^{-1}\ot p_2^*\LL^{-1}$ on $A\times A$.
The restriction $\La(\LL)|_{A\times K}$ is trivial, hence, $\La(\LL)$ descends to
a biextension $\BB$ of $A\times A/K$. Now set
$$L=p_1^*\LL\ot\BB.$$
It is easy to see that 
this line bundle on $A\times A/K$ has a structure of a $1$-cocycle with respect to the action
of $A$ on $A/K$. 
Hence, we are in the situation of Proposition \ref{cocycle-prop} with $H=A$ and $X=A/K$.
Taking $x=0\in A/K$ we obtain a central extension $G$ of $K$ by $\G_m$ with the
underlying $\G_m$-torsor $\LL|_K$ and the
right $G$-action on $\pi_*(\LL)$, where $\pi:A\to A/K$ is the projection.
One can check that if $K=K(\LL)=\ker(\phi_\LL)$ then $G\simeq G(L)^{op}$, the opposite
of the Mumford group of $L$ (see \cite[Sec.\ 12.2]{P-ab}). This explains the
seeming discrepancy between \eqref{pi(x)-pull-back-eq} and \cite[Prop.\ 12.7]{P-ab}.
\end{ex}

\subsection{Index of a symmetric isogeny}
\label{index-sec}

Recall that with every line bundle $L$ on an abelian variety $A$ one associates a symmetric
homomorphism 
$$\phi_L:A\to\hat{A}: a\mapsto t_a^*L\ot L^{-1}\in\Pic^0(A).$$
Furthermore, this induces an isomorphism
\begin{equation}\label{NS-symhom-eq}
NS(A)_\Q\to\Hom(A,\hat{A})^+_\Q
\end{equation}
of the Neron-Severi group with the space of symmetric homomorphisms tensored with $\Q$.
In the case when $\phi_L$ is an isogeny there exists a unique integer $i(L)$, $0\le i(L)\le g=\dim A$,
such that
$$H^i(A,L)=0 \text{ for } i\neq i(L),\  H^{i(L)}(A,L)\neq 0.$$
It is called the {\it index of} $L$ and can be computed as the number of positive roots of the
polynomial $P(n)=\chi(L\ot L_0^n)$, where $L_0$ is an ample line bundle on $A$ 
(see \cite[III.16]{Mum}). The index has the property
$$i(L^n)=i(L) \text{ for }\ n>0.$$
We say that an element $\phi\in\Hom(A,B)_{\Q}$ is an {\it isogeny} if it is invertible in $\Ab_{\Q}$,
or equivalently if $n\phi$ is a usual isogeny in $\Hom(A,B)$ for some $n>0$.
Using the isomorphism \eqref{NS-symhom-eq}, we can define the index of
a symmetric isogeny $\phi\in\Hom(A,\hat{A})^+_\Q$ by choosing $n>0$ such that
$n\phi=\phi_L$ and setting $i(\phi)=i(L)$. For any isomorphism $\psi:B\to A$ in
$\Ab_{\Q}$ we have
$$i(\hat{\psi}\phi\psi)=i(\phi).$$

\subsection{Orbi-abelian varieties}

We will use a natural extension of duality of abelian varieties to some length-$2$ complexes
of commutative group schemes. The
setup is similar to the duality of $1$-motives considered in \cite{De, Lau} and can be generalized
to include them (see \cite[Sec.\ 4.3]{P-FM}).

Let $\GG^{pr}_k$ denote the abelian category of commutative proper group schemes over $k$
(in characteristic zero every such group is a product of an abelian variety and a finite group).
There exists an exact duality functor
$$\Du:D^b(\GG^{pr}_k)^{op}\to D^b(\GG^{pr}_k)$$
such that $\Du(A)=\hat{A}$ for an abelian variety $A$, $\Du(G)\simeq G^*[1]$ for
a finite group scheme $G$, and $\Du^2=\Id$ (see \cite[Thm.\ 4.1.4]{P-FM}).

An {\it orbi-abelian variety} is an object $K\in D^b(\GG^{pr}_k)$ with $H^iK=0$ for 
$i\not\in\{-1,0\}$, such that $H^{-1}K$ is a finite group scheme.
Every orbi-abelian variety can be represented by a complex of the form $[G\to X]$ (in degrees $-1$
and $0$), where $X\in\GG^{pr}_k$ and $G$ is a finite group scheme.
With every orbi-abelian variety $K$ we can associate an abelian variety $H^0(K)_0$
and two finite groups: $H^{-1}(K)$ and $\pi_0H^0(K)$.
If $K$ is an orbi-abelian variety then $\Du(K)$ is also such and the corresponding
abelian varieties are dual, while the finite groups $H^{-1}(K)$ and $\pi_0H^0(K)$
get interchanged and dualized (see \cite[Lem.\ 4.2.1]{P-FM}).



\subsection{Kernels and functors}\label{conv-sec}

Recall that there is a principle originating from noncommutative geometry stating that ``correct" functors
between the derived categories of coherent sheaves on smooth projective varieties
are given by Fourier-Mukai functors (see \cite{Orlov-ker}, \cite{Toen} for some incarnations of
this principle). Such functors have form
$$\Phi_K: D^b(X)\to D^b(Y): F\mapsto Rp_{2*}(p_1^*F\ot^{{\mathbb L}} K),$$
where $K\in D^b(X\times Y)$, $p_1:X\times Y\to X$ and $p_2:X\times Y\to Y$ are the projections. 
We will refer to $K$ as a {\it kernel} representing the functor
$\Phi_K$. 

The composition of functors corresponds to an operation on kernels
that we call {\it convolution}. Namely, for $K\in D^b(X\times Y)$ and
$K'\in D^b(Y\times Z)$ we set
$$K\circ_Y K'=Rp_{13*}(p_{12}^*K\ot^{{\mathbb L}} p_{23}^*K')\in D^b(X\times Z),$$
where $p_{ij}$ are projections from $X\times Y\times Z$ to the products of two of the factors.
This operation is associative and we have a natural isomorphism of functors
$$\Phi_{K\circ_Y K'}\simeq \Phi_{K'}\circ\Phi_K$$
(see e.g. \cite[Sec.\ 2.1]{P-FM}).


\section{Representations of the Heisenberg groupoid}\label{Heis-sec}

\subsection{Symplectic setting for abelian varieties}\label{sym-sec}

This setting was introduced in \cite{P-sym} (see also \cite{P-ab}, 15.2).
The basic ``symplectic object" we work with is 
an abelian variety $X$ equipped with an isomorphism $\eta:X\to\hat{X}$
such that $\hat{\eta}=-\eta$. 

\begin{defi}
We will call $(X,\eta)$ as above a {\it symplectically selfdual abelian variety} ({\it ss-abelian variety}
for short). 
\end{defi}

The corresponding analogue of a symplectic form is the skew-symmetric biextension
$\EE=(\id\times\eta)^*\PP_X$ of $X\times X$ obtained by the pull-back from the Poincar\'e
line bundle $\PP_X$ on $X\times\hat{X}$ (see \cite{P-sym}).
Often we need to
assume that there exists a biextension $\BB$ on $X\times X$ such that $\EE$ is obtained from
$\BB$ by antisymmetrization:
$\EE\simeq\BB\otimes\si^*\BB^{-1}$, where $\si:X\to X$ is the permutation of factors.

\begin{defi}
The data $(X,\eta,\BB)$ as above is called an {\it enhanced symplectically selfdual abelian variety}
(or {\it ess-abelian variety} for short).
\end{defi}

The main example of an ess-abelian variety is 
$X_A=A\times\hat{A}$, where $A$ is an abelian variety, and 
$\BB=p_{14}^*\PP$ on $X\times X=A\times\hat{A}\times A\times\hat{A}$
(here $\PP=\PP_A$ is the Poincar\'e line bundle on $A\times\hat{A}$).
We will refer to this example as the {\it standard ess-abelian variety} associated with $A$.

We consider ess-abelian varieties as analogs of symplectic vector spaces. In \cite{P-sym}
we introduced a natural analog of the Heisenberg group in this context. 
First, let us consider the case $X_A=A\times\hat{A}$. In this case for every point
$x=(a,\xi)\in X_A(k)$ we have a functor
\begin{equation}\label{T-A-eq}
T^A_x=(\ot\PP_{\xi})\circ t_a^*: D^b(A)\to D^b(A),
\end{equation}
where $t_a:A\to A$ is the translation by $a\in A(k)$ and $\PP_{\xi}=\PP|_{A\times\{\xi\}}\in\Pic^0(A)$
is the line bundle on $A$ corresponding to $\xi\in\hat{A}(k)$.
These functors compose according to the canonical isomorphisms
$$T^A_x\circ T^A_{x'}\simeq \BB_{x,x'}\ot T^A_{x+x'}.$$
To express the algebraic dependence of these functors on $x$ we can consider similar functors
on $D^b(A\times S)$ depending on $x\in X_A(S)$, 
where $S$ is a $k$-scheme, and observe that the above isomorphisms
still hold (and are compatible with pull-backs with respect to base changes $S'\to S$).
This motivates the following definition (see \cite[Def.\ 2.1]{P-sym}).

\begin{defi} Let $(X,\BB)$ be an ess-abelian variety. The {\it Heisenberg groupoid} $\bH=\bH(X,\BB)$
associated with $(X,\BB)$ is the stack $B\G_m\times X$ equipped with the following
structure of a stack of Picard groupoids: for a $k$-scheme $S$ define the functor
$+:\bH(S)\times \bH(S)\to\bH(S)$ by 
$$(L,x)+(L',x')=(L\ot L'\ot\BB_{x,x'}, x+x'),$$
where $L$ and $L'$ are line bundles over $S$ and $x,x'\in X(S)$. 
\end{defi}

We intend to consider actions of $\bH$ on fibered categories by generalizing the example
of the action of $\bH(X_A)$ on the fibered category $S\mapsto D^b(A\times S)$ given
by the functors $T_x^A$ above (such that $L\in B\G_m(S)$ acts by tensoring with the pull-back of
$L$ to $A\times S$). 

\subsection{Isotropic and Lagrangian pairs}\label{isotr-sec}

Similarly to the case of the classical Heisenberg group associated with a symplectic vector space
we need Lagrangian subvarieties $Y$ in $X$ to define representations of $\bH$. 
Below we generalize the setup of \cite[Sec.\ 2]{P-sym} slightly in that we allow the map $Y\to X$
to have finite kernel. 

\begin{defi} A {\it generalized isotropic pair} $(Y,\a)$ consists of 
a commutative group scheme $Y$ equipped with a homomorphism $i:Y\to X$,
such that $\ker(i)$ is finite,
and of a line bundle $\a$ over $Y$ together with an isomorphism
\begin{equation}\label{isotr-eq}
\La(\a)\simeq(i\times i)^*\BB
\end{equation}
of biextensions on $Y\times Y$ inducing a cube structure on $\a$. 
In the case when $\ker(i)$ is trivial we say that $(Y,\a)$ is an {\it isotropic pair} (see \cite{P-sym}).
\end{defi}

The reason for this definition becomes clear in connection with the Heisenberg groupoid.
Namely, the data $(Y,\a)$ as above determines a homofunctor of stacks of Picard groupoids
$$Y\to\bH: y\mapsto(\a_y,i(y)).$$

Note that if $(Y,\a)$ is a generalized isotropic pair then 
the biextension $(i\times i)^*\EE$ is trivialized.
Therefore, we can view $(i\times \id)^*\EE$ as a biextension of $Y\times [Y\to X]$.
By \cite[Prop.\ 4.3.2]{P-FM},  
this biextension gives a morphism
\begin{equation}\label{lagmor}
Y\to \Du[Y\to X]
\end{equation}
that fits into a morphism of exact triangles in the derived category $D^b(\GG^{pr})$
\begin{diagram}
Y &\rTo{i} &X&\rTo{}& [Y\to X] &\rTo{}& Y[1]\\
\dTo{}&&\dTo{\eta}&&\dTo{}&&\dTo{}\\
\Du[Y\to X]&\rTo{}&\hat{X}&\rTo{\Du(i)}& \Du(Y)&\rTo{}&\Du[Y\to X][1]
\end{diagram}
Conversely, if $i:Y\to X$ is such that the composition
\begin{equation}\label{isotropic-compo-eq}
Y\rTo{i} X\rTo{\eta}\hat{X}\rTo{\Du(i)}\Du(Y)
\end{equation}
is zero then $Y$ extends to a generalized isotropic pair.
 
\begin{defi} A {\it generalized Lagrangian pair} $(Y,a)$ is a generalized isotropic pair such that
\eqref{lagmor} is an isomorphism in $D^b(\GG^{pr})$.
If $\ker(i)=0$ and $Y$ is connected then we say that $(Y,\a)$ is a Lagrangian pair and $Y$ is Lagrangian.
\end{defi}


\begin{lem}\label{lag-lem} (i)
If $(Y,\a)$ is a generalized isotropic pair then $\dim Y\le\dim X/2$. If $(Y,\a)$ is generalized Lagrangian
then $\dim Y=\dim X/2$.

\noindent
(ii) Let $Y\sub X$ be an isotropic abelian subvariety such that $\dim Y=\dim X/2$.
Then $Y$ is Lagrangian.

\noindent
(iii) Let $(Y,\a)$ be a generalized Lagrangian pair. Then $Y$ is connected if and only if $\ker(i)=0$ if
and only if $Y$ is Lagrangian.
\end{lem}

\Pf . (i) The first assertion follows from the fact that the composition
$$Y_0\rTo{i'} X\rTo{\eta} \hat{X}\rTo{\hat{i'}} \widehat{Y_0}$$
is zero and that $i'=i|_{Y_0}$ has finite kernel. 
If $(Y,\a)$ is generalized Lagrangian then we obtain from the isomorphism \eqref{lagmor}
the equality
$$\dim Y=\dim H^0\Du[Y\to X]=\dim H^0[Y\to X]=\dim X-\dim Y.$$

\noindent 
(ii)
Since $\eta(Y)$ and $\widehat{X/Y}$
are both abelian subvarieties in $\hat{X}$, it follows that $\eta(Y)=\widehat{X/Y}$ if and only if $\dim Y=\dim X/Y$.

\noindent
(iii) This follows from the duality between the finite groups
$H^{-1}[Y\to X]=\ker(i)$ and $\pi_0H^0\Du[Y\to X]\simeq\pi_0(Y)$
(see \cite[Lem.\ 4.2.1]{P-FM}).
\ed

Part (ii) of the above lemma shows that to determine whether an abelian subvariety of 
$X$ is Lagrangian we can work in the category of abelian varieties up to isogeny.

\begin{exs}\label{Lag-ex}
1. Consider the standard ess-abelian variety $X_A=A\times\hat{A}$.
Then an abelian subvariety $Y\sub X_A$, finite and surjective over $A$, gives rise
to a morphism from $A$ to $\hat{A}$ in the category $\Ab_\Q$ of abelian varieties
up to an isogeny, namely, the morphism $f_Y=p_{\hat{A}}\circ p_A^{-1}$, where
$p_A:Y\to A$ and $p_{\hat{A}}:Y\to A$ are projections.
It is easy to see that the isotropy condition on $Y$ is equivalent to the condition that
$f_Y$ is symmetric, i.e., $\hat{f}_Y=f_Y$. By Lemma \ref{lag-lem}(ii), this is equivalent
to $Y$ being Lagrangian. Conversely, with every symmetric morphism $f\in\Hom(A\to \hat{A})_{\Q}$
we can associate its graph (which will be Lagrangian)
$$\Ga(f):=\{(Na,Nf(a)) \ |\ \ a\in A\}\sub A\times\hat{A},$$
where $N>0$ is an integer such that $Nf\in\Hom(A,\hat{A})$. 
In this way we get
a one-to-one-correspondence between the set of Lagrangian subvarieties in $X_A$, finite over $A$,
and the subspace $\Hom(A,\hat{A})^+_{\Q}\sub\Hom(A,\hat{A})_{\Q}$ of symmetric homomorphisms.

\noindent
2. If $Y\sub X$ is a Lagrangian abelian subvariety then we can always choose a line
bundle $\a$ such that $(L,\a)$ is a Lagrangian pair. Indeed, the biextension 
$\BB|_{Y\times Y}$ of $Y\times Y$ is symmetric, hence, it is of the form $\La(\a)$ for some line
bundle $\a$ on $Y$ (see e.g. \cite[Thm.\ 13.7]{P-ab}).

\noindent
3. Let $X=0$. Then a generalized Lagrangian $(Y,\a)$ is a finite commutative group scheme $Y$
together with a central extension of $Y$ by $\G_m$ (given by $\a$), such that the corresponding
commutator pairing $Y\times Y\to\G_m$ induces an isomorphism $Y\to Y^*$. In other words,
this gives a finite Heisenberg group scheme. 
\end{exs}

\begin{lem}\label{conn-comp-lem}
If $(Y,\a)$ is a generalized isotropic pair such that $\dim Y=\dim X/2$ then
$i(Y_0)\sub X$ is a Lagrangian abelian subvariety, where $Y_0\sub Y$ is the connected component
of $0$ in $Y$.
\end{lem}

\Pf . By definition, the composition $Y\to X\wt{\to} \hat{X}\to \Du(Y)$ is zero. Hence,
the same is true for the composition 
$$Y_0\to X\wt{\to}\hat{X}\to\widehat{Y_0}.$$
Since the natural morphisms $Y_0\to i(Y_0)$ and $\widehat{i(Y_0)}\to\widehat{Y_0}$
are isogenies, we derive that the composition 
$$i(Y_0)\to X\wt{\to}\hat{X}\to\widehat{i(Y_0)}$$
is zero. By Lemma \ref{lag-lem}(ii), this implies that $i(Y_0)$ is Lagrangian.
\ed


\begin{defi} Let us say that generalized Lagrangian pairs $(Y,\a)$ and $(Z,\b)$ are {\it transversal}
if the images of $Y$ and $Z$ generate $X$, or equivalently if $Y\times_X Z$ is finite.
\end{defi}

\begin{lem}\label{transv-lem} 
(i) For any pair of Lagrangian subvarieties $Y$ and $Z$ in $X$ there exists a Lagrangian subvariety
$Q\sub X$ transversal to both $Y$ and $Z$. 

\noindent
(ii) Now assume in addition that $X=X_A$ (equipped with the standard symplectic biextension) 
and $Y=\{0\}\times\hat{A}$. Let $\LL$ be an ample line bundle  on $\hat{A}$.
Then for almost all $n\in\Z$ (i.e., for all except for a finite
number) the Lagrangian subvariety $\Ga(n\phi_\LL)\sub \hat{A}\times A=X_A$ 
(see Example \ref{Lag-ex}.1) is transversal to $Z$.
\end{lem} 

\Pf . (i) To prove the first assertion we can argue in the semisimple category $\Ab_{\Q}$
of abelian varieties up to isogeny (since we have natural analogs of the relevant notions
in $\Ab_{\Q}$). In particular, in this category 
the condition of transversality of $Q$ to $Y$ and $Z$ becomes simply 
$Q\cap Y=Q\cap Z=0$.
We can assume that $X\simeq Y\oplus \hat{Y}$, where
both $Y$ and  $\hat{Y}$ are Lagrangian in $X$ and the self-duality of $X$ restricts
to the duality between $Y$ and $\hat{Y}$ (see \cite[Lemma 1.1]{P-sym}).  
Let $T=Y\cap Z$. We can write 
\begin{equation}\label{Y-T-P-dec-eq}
Y=T\oplus P, \ \ X=T\oplus P\oplus \hat{T}\oplus\hat{P},
\end{equation}
so that 
\begin{equation}\label{Z-T-P-eq}
Z=\{(t,\phi(\xi),0,\xi) |\ t\in T, \xi\in\hat{P}\},
\end{equation}
where $\phi:\hat{P}\to P$ is a symmetric morphism.
Then choosing a symmetric isomorphism $\psi:\hat{P}\to P$ in $\Ab_{\Q}$
we can set
$$Q=\{(0,(\phi+\psi)(\xi),\eta,\xi) |\ \eta\in\hat{T}, \xi\in\hat{P}\}.$$

\noindent
(ii) 
As in part (i) we consider the decompositions \eqref{Y-T-P-dec-eq} in $\Ab_{\Q}$ (with $Y=\hat{A}$)
so that
\eqref{Z-T-P-eq} holds (still in $\Ab_{\Q}$). 
Let us write the homomorphism 
$\phi_\LL:Y\to \hat{Y}$ in the form
$$\phi_\LL=\left(\begin{matrix} \a & \b \\ \ga & \de\end{matrix}\right),$$
with $\a\in\Hom(T,\hat{T})_\Q$, $\b\in\Hom(P,\hat{T})_\Q$, $\ga\in\Hom(T,\hat{P})_\Q$ and
$\de\in\Hom(P,\hat{P})_\Q$. It is easy to see that the condition of transversality of $\Ga(n\phi_\LL)$
with $Z$ is equivalent to the condition that $n\de-\phi$ is an isogeny (i.e., invertible in $\Ab_\Q$).
Note that some multiple of $\de$ is the symmetric homomorphism $P\to\hat{P}$ 
corresponding to the restriction of $\LL$ to $P$. Since this restriction is still ample, it follows
that $\de$ is an isogeny. Thus, $\deg(n\de-\phi)$ is a polynomial in $n$ with nonzero top degree coefficient
$\deg(\de)$, so $\deg(n\de-\phi)\neq 0$ for almost all $n\in\Z$.
\ed

The construction of the following proposition is a straightforward
generalization of \cite[Lem.\ 15.4]{P-ab}
(see also \cite[Sec.\ 3]{P-sym}).
Recall that given a biextension $\LL$ of $Y\times Z$ by $\G_m$ and
a pair of homomorphisms $A\to Y$, $B\to Z$, together with 
trivializations of the pull-backs of $\LL$ to $A\times Z$ and $Y\times B$,
we obtain a pairing 
$$A\times B\to\G_m$$
measuring the difference between two induced trivializations of the pull-back
of $\LL$ to $A\times B$ (cf. \cite[10.4]{P-ab}). 

\begin{prop}\label{lag-int-prop}
Let $(i:Y\to X,\a)$ and $(j:Z\to X,\b)$ be generalized Lagrangian pairs such that
$Y\times_X Z$ is finite.
Then there is a natural central extension $G$ of $K=Y\times_X Z$ by $\G_m$
with the underlying line bundle $\b_z\ot \a_y^{-1}$ over $K$.
Moreover, $G$ is a Heisenberg group scheme and the corresponding commutator
form $K\times K\to\G_m$ is the pairing associated with the biextension
$(i\times j)^*\EE$ of $Y\times Z$ and with natural trivializations of its pull-backs
to $K\times Z$ and $Y\times K$.
\end{prop}

\Pf . Pulling back the isomorphisms
\begin{equation}\label{La-i-j-eq}
\La(\a)\simeq(i\times i)^*\BB, \ \ \La(\b)\simeq(j\times j)^*\BB
\end{equation}
(see \eqref{isotr-eq}) 
to $K\times K$ we obtain an isomorphism $\La(\b_z\ot \a_y^{-1})\simeq\OO_{K\times K}$
that gives a central extension structure on the corresponding $\G_m$-torsor $G\to K$.
The isomorphisms \eqref{La-i-j-eq} also induce
trivializations of the pull-backs of $\EE$ to $Y\times Y$ and to $Z\times Z$. 
Hence, we have natural trivializations of the pull-backs
of $\EE$ to $K\times Z$ and $Y\times K$ whose difference gives a pairing 
$$e:K\times K\to\G_m.$$
Note that $\EE$ induces a duality between $Y$ and $[Y\to X]\simeq [K\to Z]$ and
between $Z$ and $[Z\to X]\simeq [K\to Y]$. Thus, the dual of the exact triangle
$$K\to Z\to [K\to Z]\to\ldots$$
is the exact triangle
$$Y\to [K\to Y]\to K[1]\to\ldots$$
and we obtain an isomorphism $K\simeq K^*$. It is easy to see that it is given by the pairing
$e$. Furthermore, the same reasoning as in \cite[Lem.\ 15.4]{P-ab} shows that $e$ is
the commutator form of the central extension $G\to K$.
\ed

\subsection{Representations associated with Lagrangian pairs and intertwining functors}

Let $(Y,\a)$ be an isotropic pair. Then the data 
$(Y\to X,\BB|_{Y\times X},\a)$ is similar to the twisting data considered in \cite[Sec.\ 3.7]{P-FM}
with the only difference that $Y$ is not a finite group scheme. In particular,
we still have the $1$-cocycle of $Y$ with values in $\PPic(X)$ defined
by
\begin{equation}\label{Y-cocycle-eq}
\LL=p_1^*\a^{-1}\otimes\BB|_{Y\times X}^{-1}.
\end{equation}

Recall that in \cite[Sec.\ 2]{P-sym} we defined 
the category $\FF(Y,\a)$ as the category of $\LL$-twisted $Y$-equivariant objects
in $D^b(X)$ (see Sec.\ \ref{twisted-sec}). 
Explicitly, the objects of $\FF(Y,\a)$ are objects $F\in D^b(X)$ equipped with isomorphisms
$$F_{x+y}\simeq\a_y^{-1}\ot\BB_{y,x}^{-1}\ot F_x$$
on $X\times Y$, satisfying the natural cocycle condition on $X\times Y\times Y$.

For every $k$-scheme $S$ we can define similarly the category 
$\FF(Y,\a)_S$ as the category of $\LL$-twisted $Y$-equivariant objects in $D^b(X\times S)$,
so that $S\mapsto \FF(Y,\a)_S$ is a fibered category with respect to the natural
pull-back functors.
The Heisenberg groupoid $\bH$ naturally acts on this fibered category. Namely,
for every $S$ we have an action of $\bH(S)$ on $\FF(Y,\a)_S$ given by the functors
$$T_{L,x}(F)_{x'}=L\ot\BB_{x',x}\ot F_{x+x'},$$
and these actions are compatible with the pull-backs with respect to morphisms $S\to S'$.
We also set $T_x=T_{\OO_S,x}$ for $x\in X(S)$.

\begin{rems}\label{F-Y-rem} 
1.
In the case when the projection $X\to X/Y$ has a section $s:X/Y\to X$ there is a natural
equivalence $\FF(Y,\a)\simeq D^b(X/Y)$ obtained by restricting $F\in\FF(Y,\a)$ to 
$s(X/Y)\sub X$. Furthermore, the functors $T_{L,x}$ are given by some kernels on
$X/Y\times X/Y$, so we get a {\it kernel representation} in the sense of \cite[Sec.\ 3.2]{P-FM}
(see also Example \ref{main-rep-ex} below).

\noindent
2. The above definition of $\FF(Y,\a)$ is not quite adequate in general. 
For example, it is not clear whether $\FF(Y,\a)$ has a triangulated structure such that
the embedding $\FF(Y,\a)\to D^b(X)$ is an exact functor.
If $\cha(k)=0$ then $\FF(Y,\a)$ can be identified with the (bounded) derived category of
modules over some Azumaya algebra over $X/Y$ (see \cite[Prop.\ 3.1]{P-sym}).
This identification uses the derived descent with respect to finite flat morphisms 
(see \cite[Appendix]{P-sym}\footnote{The assumption $\cha(k)=0$
was erroneously omitted in \cite{P-sym}}, \cite{El}). 
A more adequate replacement for $\FF(Y,\a)$ is the derived category of the abelian
category $\FF(Y,\a)\cap\Coh(X)=\Coh^\LL_Y(X)$ (cr.\ Remark \ref{equiv-rem} below).
However, we don't need this since we will mostly work with objects of $\FF(Y,\a)$ that are cohomologically pure. 
\end{rems}

\begin{ex}\label{main-rep-ex} Consider the standard ess-abelian variety $X_A=A\times\hat{A}$
with $\BB=p_{14}^*\PP$ and take $Y=\{0\}\times\hat{A}\sub X_A$, $\a=\OO$. Then using
the section $A\times\{0\}\sub X_A$ of the projection $X_A\to A$ we obtain an equivalence of 
the category $\FF(Y,\a)$ with $D^b(A)$. One immediately checks that under this equivalence
the functors $T_x$ giving the action of $\bH(k)$ on $\FF(Y,\a)$ become precisely the functors
$T_x^A$ defined by \eqref{T-A-eq}. 
We will denote by $\bbT^A_x\in D^b(A\times A)$ the kernel giving the functor $T_x^A$.
Explicitly, for $x=(a,\xi)\in A\times\hat{A}$,
\begin{equation}\label{bbT-eq}
\bbT^A_x=(t_a,\id_A)_*\PP_{\xi}.
\end{equation}
\end{ex}

\begin{rem}\label{moduli-rem}
Given an object $F\in D^b(A)$ we can consider the subset $Y_F$ of
$x\in X_A$ such that $T_x^A(F)\simeq F$. 
Assuming that $F$ is endosimple and $\Ext^i(F,F)=0$ for $i<0$, one can check
that $Y_F$ is an algebraic subgroup in $X_A$ and that
we have a line bundle $\a$ on $Y_F$ and an isomorphism
$$T_y^A(F)\simeq \a_y\ot F$$
over $Y_F\times A$ (this can be done using results of \cite{Inaba}; the case of
vector bundles is considered in \cite[Sec.\ 1]{Mukai}).
This easily implies that $(Y_F,\a)$ is an isotropic pair. In particular, $\dim Y_F\le \dim A$.
\end{rem}


As in the classical case of the Heisenberg group, one expects that the representations $\FF(Y,\a)$ associated with different Lagrangian pairs $(Y,\a)$ are equivalent. 
This can be proved under some additional assumptions (due to the need to use the derived
descent, see Remark \ref{F-Y-rem}.2).

\begin{thm}\label{main-eq-thm} 
Let $(Y,\a)$ and $(Z,\b)$ be Lagrangian pairs for an ess-abelian variety $X$.
Assume that either $\cha(k)=0$ or $Y\cap Z=0$.
Then there is an equivalence of fibered categories
$\FF(Y,\a)_S\simeq \FF(Z,\b)_S$ compatible with the $\bH$-action. 
\end{thm}

\Pf . In both cases the proof of \cite[Thm.\ 4.3]{P-sym} works. 
The assumption $\cha(k)=0$ allows one to use the derived descent, while in the case
$Y\cap Z=0$ the descent problem does not arise.
\ed

\begin{rem}\label{equiv-rem}
In this paper we use only categories $\FF(Y,\a)$ associated with Lagrangian pairs. However,
the above equivalence can also be constructed for generalized Lagrangian pairs
and the extra assumptions in Theorem \ref{main-eq-thm} can be removed
if we replace each category $\FF(Y,\a)$ with the derived
category of $\FF(Y,\a)\cap\Coh(X)$ and use Theorem 4.5.1 of \cite{P-FM}. 
\end{rem}

%

\subsection{Invariants of a generalized Lagrangian pair}\label{inv-sec}

Consider the representation of $\bH$ on the fibered category $S\mapsto \FF(Y,\a)_S$,
associated with a Lagrangian pair $(Y,\a)$.
Now let $(j:Z\to X, \b)$ be a generalized isotropic pair.
Then we have a homofunctor 
$$\bj:Z\to\bH: z\mapsto (\b_z, j(z))$$
of Picard stacks, 
so it makes sense to consider the
category $\FF(Y,\a)^{(Z,\b)}$ of $(Z,\b)$-invariants in $\FF(Y,\a)$.
By definition, the objects of this category are objects $F\in\FF(Y,\a)$ equipped
with a collection of isomorphisms
$$T_{\bj(z)}(F_S)\rTo{\phi(z)} F_S,$$
for all $k$-schemes $S$ and all $z\in Z(S)$, where $F_S\in\FF(Y,\a)_S$ is the
pull-back of $F$, compatible with pull-backs with respect to base changes $S\to S'$
and satisfying the following cocycle condition
\begin{diagram}
T_{\bj(z_1)+\bj(z_2)}(F_S) &\rTo{T_{\bj(z_1)}\phi(z_2)}& T_{\bj(z_1)}(F_S)\\
\dTo{}&&\dTo{\phi(z_1)}\\
T_{\bj(z_1+z_2)}(F_S)&\rTo{\phi(z_1+z_2)}& F_S
\end{diagram}
One can easily check that the collection $(\phi(z))$ is determined by the single element
$\phi(z^{un})$, where $z^{un}\in Z(Z)$ is the tautological $Z$-point of $Z$. 
Thus, $\FF(Y,\a)^{(Z,\b)}$ can be also described as
the category of objects $F\in\FF(Y,\a)$ equipped
with isomorphisms
$$F_{j(z)+x}\simeq\b_z^{-1}\ot\BB_{x,j(z)}^{-1}\ot F_x,$$
over $Z\times X$, where $z\in Z$, $x\in X$, satisfying the natural cocycle condition 
over $Z\times Z\times X$ and the following compatibility for $F_{j(z)+x+y}$, where $y\in Y$, $z\in Z$:
\begin{diagram}
F_{j(z)+x+y}&\rTo{} &\b_z^{-1}\ot \BB_{x+y,j(z)}^{-1}\ot F_{x+y}\\
\dTo{}&&\dTo{}\\
\a_y^{-1}\ot\BB_{y,j(z)+x}^{-1}\ot F_{j(z)+x}&\rTo{}&\a_y^{-1}\ot\b_z^{-1}\ot\BB_{y,j(z)}^{-1}\ot\BB_{y,x}^{-1}\ot
\BB_{x,j(z)}^{-1}\ot F_x
\end{diagram}
where the horizontal (resp., vertical) isomorphisms correspond to $(Z,\b)$-invariance 
(resp., come from the definition of $\FF(Y,\a)$).


\begin{ex}\label{homog-ex} 
As in Example \ref{main-rep-ex} let us consider the standard ess-abelian variety
$X=X_A$ and $Y=\{0\}\times\hat{A}\sub X_A$, so that $\FF(Y,\OO)\simeq D^b(A)$.
Then the skyscraper sheaf $k(0)\in D^b(A)$ is invariant with respect to 
$(Y,\OO)$. Now let $F$ be a nontrivial extension of $k(0)$ by itself, so $F$ is isomorphic
to the structure sheaf of a length $2$ subscheme $S\sub A$.
Then the fact that the restriction of the Poincar\'e bundle to $S\times\hat{A}$ is nontrivial implies
that $F$ is not $(Y,\OO)$-invariant, even though $F\ot\PP_{\xi}\simeq F$ for every
$\xi\in\hat{A}(\xi)$. The Fourier dual of this example is a bundle $V$ which is a nontrivial extension
of $\OO_A$ by itself. The bundle $V$ is homogeneous, i.e., $t_a^*V\simeq V$ for every
$a\in A$. However, it is not invariant with respect to $A\times\{0\}\sub X_A$ (and $\b=\OO$)
according to our definition.
\end{ex}

By Theorem \ref{main-eq-thm}, if $(Y',\a')$ is another Lagrangian pair then
the categories $\FF(Y,\a)^{Z,\b}$ and $\FF(Y',\a')^{Z,\b}$ are equivalent, 
provided either $\cha(k)=0$ or $Y\cap Y'=0$.
This observation will often allow us to reduce the study of $\FF(Y,a)^{Z,\b}$ to
the case when $Y$ and $Z$ are transversal.

\begin{lem}\label{beta-change-lem}
Let $\b'=\b\ot \EE|_{Z\times\{x\}}$ for some $x\in X$
then the functor $T_x:\FF(Y,\a)\to\FF(Y,\a)$ induces an equivalence
$$t(Z,x):\FF(Y,\a)^{(Z,\b')}\to\FF(Y,\a)^{(Z,\b)}.$$
\end{lem}

\Pf . Let $\bj':Z\to\bH: z\mapsto(\b'_z, j(z))$ be the homofunctor associated with $\b'$.
Then we have 
$$T_{\bj'(z)}=\EE_{j(z),x}\ot T_{\bj(z)}.$$
Hence, the structure of $(Z,\b')$-invariance for $F\in\FF(Y,\a)$ induces isomorphisms
$$T_{\bj(z)}T_x(F)\simeq \EE_{j(z),x}\ot T_xT_{\bj(z)}(F)\simeq
T_x T_{\bj'(z)}(F)\to T_x(F),$$
i.e., the structure of $(Z,\b)$-invariance for $T_x(F)$.
\ed

We are going to show that under some technical assumptions the category 
$\FF(Y,\a)^{(Z,\b)}$ is generated by a single coherent sheaf in the sense of the following definition.

\begin{defi} Let $\TT$ be a triangulated category, $\FF\sub\TT$ an additive (not necessarily
full) subcategory closed under shifts $X\mapsto X[i]$, $i\in\Z$. 
For objects $F,F_1,F_2\in\FF$
we say that $F$ is an {\it $\FF$-extension} of $F_2$ by $F_1$ if there exists
a triangle  
$$F_1\to F\to F_2\to F_1[1]$$
in $\FF$ which is an exact triangle in $\TT$. We say that
$F\in\FF$ has an {\it $\FF$-filtration with successive quotients} $Q_1,\ldots,Q_n\in\FF$
if there exists a collection of objects $(F_{ij})$, $0\le i<j\le n$, in $\FF$, such that 
$F_{i-1,i}=Q_i$ for $i=1,\ldots,n$, $F_{0n}=F$, and for any $i<j<k$, 
$F_{ik}$ is an $\FF$-extension of $F_{jk}$ by $F_{ij}$.
We say that an object
$F_0\in\FF$ is a {\it $t$-generator of $\FF$ with respect to $\TT$} if 
$\Hom_{\TT}(F_0,F_0[i])=0$ for $i<0$,
and every nonzero object of $\FF$ admits an $\FF$-filtration with successive quotients
$F_0^{\oplus m_1}[k_1],\ldots,F_0^{\oplus m_n}[k_n]$ for some $n>0$, where
$k_1>k_2>\ldots>k_n$.
\end{defi}

Note that the above definition of a $t$-generator mimics the situation when $\FF$ has
a $t$-structure with the heart consisting of direct sums of $F_0$, however, it uses
the triangulated structure of an ambient category $\TT$ (since $\FF$ does not have to
be triangulated).

\begin{lem}\label{t-gen-lem} 
Let $F_0$ be a $t$-generator of $\FF$ with respect to $\TT$. 

\noindent
(i) Suppose
$A\in\FF$ (resp., $B\in\FF$) has an $\FF$-filtration with successive quotients
$$F_0^{\oplus m_1}[a_1],\ldots,F_0^{\oplus m_n}[a_n] \ \text{ (resp., }\ 
F_0^{\oplus m'_1}[b_1],\ldots,F_0^{\oplus m'_p}[b_p] \text{ )}$$ 
such that
$a_1>\ldots>a_n$ (resp., $b_1>\ldots>b_p$).
Then
$$\Hom_{\FF}(A,B[a_n-b_1])\neq 0.$$

\noindent
(ii) Let $F\in\FF$ be an object such that $\Hom_{\FF}(F,F[i])=0$ for $i<0$. Then 
$F\simeq F_0^{\oplus m}$ in $\FF$.

\noindent
(iii) Suppose
$F_0$ is an $\FF$-extension of $B\in\FF$ by $A\in\FF$, where 
$\Hom_{\TT}(A,B[i])=0$ for $i\le 0$. Then either $A=0$ or $B=0$.
\end{lem}

\Pf . (i) By definition, we have an exact triangle
$$A'\to A\rTo{\pi} F_0^{\oplus m_n}[a_n]\to A'[1]$$
where $A'$ has an $\FF$-filtration with successive quotients of the form $F_0^{\oplus m}[i]$
such that $i>a_n$. Hence, $\Hom_{\TT}(A'[1],F_0[a_n])=0$ and the map
$$\Hom_{\FF}(F_0^{\oplus m_n}[a_n],F_0[a_n])\to\Hom_{\FF}(A,F_0[a_n])$$
induced by $\pi$ is injective. In particular, $\Hom_{\FF}(A,F_0[a_n])\neq 0$.
Next, let us consider an exact triangle
$$F_0^{\oplus m'_1}[a_n]\rTo{\iota} B[a_n-b_1]\to B'\to\ldots$$
where $B'$ has an $\FF$-filtration with successive quotients of the form $F_0^{\oplus m}[i]$,
such that $i<a_n$. It follows that $\Hom_{\TT}(A,B'[-1])=0$, so the map 
$$\Hom_{\FF}(A,F_0^{\oplus m'_1}[a_n])\to \Hom_{\FF}(A,B[a_n-b_1])$$
induced by $\iota$ is injective. Since the source space is non-zero, this implies the result.

\noindent
(ii) We know that $F$ has an $\FF$-filtration with successive quotients
$F_0^{\oplus m_1}[k_1],\ldots,F_0^{\oplus m_n}[k_n]$ such that
$k_1>\ldots>k_n$. Applying (i) to $A=B=F$ we deduce that $k_n\ge k_1$, so in
fact, $n=1$.

\noindent
(iii) Assume that $A$ and $B$ have $\FF$-filtrations as in (i). Then the result of (i) together
with our vanishing assumption implies that $b_1<a_n$. On the other hand, in
order to have a nonzero morphism $A\to F_0$ (resp., $F_0\to B$) we should have
$a_n\le 0$ (resp., $b_1\ge 0$). Thus, either $\Hom_{\TT}(A,F_0)=0$ or
$\Hom_{\TT}(F_0,B)=0$. Assume that $\Hom_{\TT}(A,F_0)=0$ (the second case is considered
similarly). Then in $\TT$ one has an isomorphism $B\simeq F_0\oplus A[1]$, so we obtain that
$$\Hom_{\TT}(A,B[-1])=\Hom_{\TT}(A,A)\neq 0,$$
a contradiction.
\ed

\begin{thm}\label{lag-inv-thm} Let $(Y,\a)$ be a Lagrangian pair for an ess-abelian variety $X$, and
let $(j:Z\to X,\b)$ be a generalized Lagrangian pair for $X$. Assume that the following condition
is satisfied:

\noindent
{\bf $(\star)$
either $Y$ and $Z$ are transversal, or $X=X_A$ and $Y=\{0\}\times\hat{A}$, or $\cha(k)=0$}.

\noindent
Then there exists a coherent sheaf
$$S_{Z,\b}=S_{Z,\b}(Y,\a)\in \FF(Y,\a)^{(Z,\b)}\cap\Coh(X)$$ 
which is a $t$-generator of $\FF(Y,\a)^{(Z,\b)}$ with respect to $\FF(Y,\a)$. Furthermore,
$S_{Z,\b}$ is an  endosimple object of $\FF(Y,\a)^{(Z,\b)}$.
Every coherent sheaf in $\FF(Y,\a)^{(Z,\b)}\cap\Coh(X)$ is a direct sum of several
copies of $S_{Z,\b}$.
In the case when
$Y$ and $Z$ are transversal we have
\begin{equation}\label{lag-inv-formula}
S_{Z,\b}(Y,\a)=(s_{Y,Z})_*(\a^{-1}\boxtimes\b^{-1}\otimes (i\times j)^*\BB^{-1})^I,
\end{equation}
where $i:Y\to X$ is the embedding map, $s_{Y,Z}:Y\times Z\to X$ is the natural map and
$I$ is a maximal isotropic subgroup in $\ker(s_{Y,Z})\simeq Y\times_X Z$
with respect to the commutator pairing from Proposition \ref{lag-int-prop}.
We also have
\begin{equation}\label{End-F-Y-eq}
\dim_k\End_{\FF(Y,\a)}(S_{Z,\b})=|\pi_0(Z)|.
\end{equation}
\end{thm}

\Pf . Assume first that $Y$ and $Z$ are transversal and let us show the existence
of an endosimple $t$-generator in this case.
Recall that an object $F\in\FF(Y,\a)$ is equipped with an isomorphism
$$F_{x+y}\simeq\a_y^{-1}\ot \BB_{y,x}^{-1}\ot F_x,$$
where $y\in Y$, $x\in X$.
On the other hand, the condition of $(Z,\b)$-invariance is an isomorphism
$$F_{j(z)+x}\simeq\b_z^{-1}\ot \BB_{x,j(z)}^{-1}\ot F_x,$$
where $z\in Z$, $x\in X$. 
Combining these two isomorphisms we get an isomorphism
\begin{equation}\label{L-twisted-main-eq}
F_{j(z)+x+y}\simeq L_{(y,z),x}\ot F_x,
\end{equation}
over $(Y\times Z)\times X$, where 
\begin{equation}\label{L-cocycle}
L_{(y,z),x}=\a_y^{-1}\ot \b_z^{-1}\ot \BB_{x,j(z)}^{-1}\ot\BB_{y,j(z)}^{-1}\ot\BB_{y,x}^{-1}.
\end{equation}
has a natural structure of a $1$-cocycle for the action of $Y\times Z$ on $X$
induced by the homomorphism $s_{Y,Z}:Y\times Z\to X$. 
Assume in addition that $F$ is a coherent sheaf.
Then the compatibilities in the definition of $\FF(Y,\a)^{(Z,\b)}$
reduce to the condition that \eqref{L-twisted-main-eq} 
gives $F$ a structure of an $L$-twisted $Y\times Z$-equivariant
coherent sheaf on $X$. By Proposition \ref{cocycle-prop}, the category of such sheaves is equivalent
to the category of weight one representations of
the central extension $G$ of $K=\ker(s_{Y,Z})$ by $\G_m$ given by the restriction
$L|_{K\times\{0\}}$. We have an identification 
$$\iota:Y\times_X Z\to K:(y,z)\mapsto (y,-z)$$
such that 
$$\iota^*(L|_{K\times\{0\}})_{y,z}\simeq \a_y^{-1}\ot\b_{-z}^{-1}\ot \BB_{j(z),-j(z)}^{-1}\simeq
\a_y^{-1}\ot\b_{j(z)},$$
where the last isomorphism comes from the isomorphism $\La(\b)\simeq (j\times j)^*\BB$.
Hence, $G$ can be identified with the Heisenberg group of Proposition \ref{lag-int-prop}.
Let us define $S_{Z,\b}\in\FF(Y,\a)^{Z,\b}\cap\Coh(X)$ as the sheaf corresponding
to the Schr\"odinger representation of $G$. Note that it is an endosimple object of $\FF(Y,\a)^{Z,\b}$.
The formula \eqref{lag-inv-formula} follows immediately from Lemma \ref{push-forward-lem}.
For any object $F\in \FF(Y,\a)^{Z,\b}\sub D^b(X)$ the triangles
\begin{equation}\label{devissage-tr-eq}
\tau_{\le n-1}F\to F\to \tau_{\ge n}F\to \tau_{\le n-1}F[1]
\end{equation}
(where $\tau_\bullet$ are the truncation functors associated with the standard $t$-structure)
automatically belong to $\FF(Y,\a)^{Z,\b}$. Hence, the standard devissage of $F$ into its cohomology
sheaves shows that $S_{Z,\b}$ is a $t$-generator of $\FF(Y,\a)^{Z,\b}$.

The existence of an endosimple $t$-generator of $\FF(Y,\a)^{Z,\b}$ with respect to $\FF(Y,\a)$
in the case when
$Y$ and $Z$ are not necessarily transversal follows now from Theorem \ref{main-eq-thm}. 
More precisely, if $\cha(k)=0$ then we can choose a Lagrangian pair $(Y',\a')$, transversal
to $Z$ (see Lemma \ref{transv-lem}(i)) 
and use an equivalence of fibered categories $\FF(Y,\a)_S\simeq\FF(Y',\a')_S$
compatible with the $\bH$-action. In the case when
$Y=\{0\}\times\hat{A}\sub X_A=X$ we can assume that $Y'\cap Y=0$ (see Lemma \ref{transv-lem}(ii)),
so we again have an equivalence $\FF(Y,\a)_S\simeq\FF(Y',\a')_S$.

Let us prove that such a $t$-generator $F_0$ of $\FF(Y,\a)^{Z,\b}$ with respect to $\FF(Y,\a)$ 
is necessarily cohomologically pure.
Suppose $F_0$ is not pure and consider the exact triangle \eqref{devissage-tr-eq}
with $A=\tau_{\le n-1}F_0\neq 0$ and $B=\tau_{\ge n}F_0\neq 0$. Then
$\Hom(A,B[i])=0$ for $i\le 0$, hence, the same space of morphisms in $\FF(Y,\a)$
also vanishes. But this gives a contradiction with Lemma \ref{t-gen-lem}(iii).

Thus, $F_0$ can be shifted to become a coherent sheaf, and we set
$S_{Z,\b}=F_0$.
The fact that every $F\in\FF(Y,\a)^{(Z,\b)}\cap\Coh(X)$ is a direct sum of several copies of $S_{Z,\b}$
follows from Lemma \ref{t-gen-lem}(ii).

To prove \eqref{End-F-Y-eq} we can again assume that $Y$ and $Z$ are transversal.
Let $V$ be a Schr\"odinger representation of the Heisenberg group scheme $G$ considered above.
By Lemma \ref{push-forward-lem}, endomorphisms of $F=S_{Z,\b}$ can be identified with
$K$-invariants in endomorphisms of 
$$s_{Y,Z}^*F\simeq V\ot L|_{Y\times Z\times\{0\}},$$
i.e., with $K$-invariants in $\End_k(V)\ot H^0(Y\times Z,\OO)$.
But $\End_k(V)$ is isomorphic to the space of functions on $K$ (see \eqref{End-eq}), so
$$\End(S_{Z,\b})\simeq H^0(K\times Y\times Z,\OO)^K\simeq H^0(K\times\pi_0(Z),\OO)^K.$$
The dimension of the latter space is $|\pi_0(Z)|$.
\ed

\begin{ex} Recall that in the case when $X$ is the standard ess-abelian variety $X_A=A\times\hat{A}$
and $Y=0\times\hat{A}$ we have an equivalence $\FF(Y,\OO)\simeq D^b(A)$.
Consider the Lagrangian subvariety $Z=L_\phi\sub A\times\hat{A}$ associated with a symmetric
homomorphism 
$\phi\in\Hom(A,\hat{A})^+_\Q$ (see Example \ref{Lag-ex}.1) and extend it in some way to a Lagrangian
pair $(Z,\b)$. Then by Theorem \ref{lag-inv-thm}, the LI-sheaf $E=S_{Z,\b}\in D^b(A)$ is an 
endosimple vector bundle.
Furthermore, this bundle is semihomogeneous in the sense of \cite{Mukai}, and
$\phi$ corresponds to the element
$\det(E)/\rk(E)\in NS(A)_\Q$ under the isomorphism \eqref{NS-symhom-eq} 
(see \cite[Prop.\ 7.7]{Mukai}).
\end{ex}

In the next proposition we study the relation between invariants of a generalized Lagrangian
$j:Z\to X$ in $\FF(Y,\a)$ and invariants of the corresponding Lagrangian $j(Z_0)\sub X$.

\begin{prop}\label{lag-inv-prop} 
Let $(Y,\a)$ be a Lagrangian pair. 

\noindent
(i) Let $(\ov{Z},\ov{\b})$ be a Lagrangian pair, transversal to $Y$,
and let $j:Z_0\to \ov{Z}$ be a surjective homomorphism 
from an abelian variety $Z_0$ with finite kernel $C$, 
so that we have the generalized isotropic pair $(Z_0,\b)$, where
$\b=j^*\ov{\b}$.
Consider the subgroup 
$$Z'=\ker(X\rTo{\eta}\hat{X}\to \hat{Z_0})\sub X.$$
Then we have a natural duality
$$\chi:Z'/\ov{Z}\rTo{\sim} C^*.$$
and an equivalence
$$\Phi:\FF(Y,\a)^{(Z_0,\b)}\cap\Coh(X)\rTo{\sim} C-\rep,$$
such that for $F\in \FF(Y,\a)^{(Z_0,\b)}$ and $z'\in Z'$ one has
\begin{equation}\label{translation-chi-isom}
\Phi (T_{z'} F)\simeq \chi(z'\mod\ov{Z})\ot \Phi(F).
\end{equation}
Furthermore, $\Phi(S_{\ov{Z},\ov{\b}})$ is isomorphic to the trivial
one-dimensional representation of $C$.

\noindent
(ii) Let $(j:Z\to X,\b)$ be a generalized Lagrangian pair and let $Z_0$ be the connected 
component of zero in $Z$. Assume that the condition $(\star)$ from Theorem \ref{lag-inv-thm}
is satisfied. Consider the finite group scheme $\Pi=\pi_0(j(Z))=j(Z)/j(Z_0)$ as a subgroup 
of $X/j(Z_0)$.
Assume first that $\Pi$ is reduced.
Then there exists a $\Pi$-coset $\Pi'\sub X/j(Z_0)$ such that 
\begin{equation}\label{S-Z-Z0-coset-eq}
S_{Z,\b}\simeq \left(\bigoplus_{\ov{x}\in\Pi'} T_{\ov{x}} S_{j(Z_0),\ov{\b}}\right)^{\oplus N_Z}
\end{equation}
in $\FF(Y,\a)$,
where $\ov{\b}$ is chosen in such a way that $(j(Z_0),\ov{\b})$ is a Lagrangian pair, and 
$$N_Z=\frac{|\pi_0(Z)|^{1/2}}{|\pi_0(j(Z))|^{1/2}}$$
($N_Z$ is always an integer).
If $\Pi$ is not reduced then there is still a $\Pi$-coset $\Pi'\sub X/j(Z_0)$
such that the composition factors of $S_{Z,\b}$ in $\FF(Y,\a)^{(Z_0,\b)}\cap\Coh(X)$
(which is a finite length abelian category)
are $(T_{\ov{x}} S_{j(Z_0),\ov{\b}})$, $\ov{x}\in\Pi'$,
each with multiplicity $N_Z\cdot \frac{|\Pi|}{\# \Pi(k)}$.
\end{prop}

\Pf . (i) As in the proof of Theorem \ref{lag-inv-thm} 
the line bundle $L$ on $(Y\times Z_0)\times X$ given by 
\eqref{L-cocycle} has a structure of a $1$-cocycle
for the (transitive) action of $Y\times Z_0$ on $X$. Furthermore, this $1$-cocycle is
a pull-back of a similar cocycle $\ov{L}$ on $(Y\times \ov{Z})\times X$. Therefore,
the central extension $G$ of $K=\ker(s_{Y,Z_0})$ by $\G_m$ is the pull-back of the similar
central extension $\ov{G}$ of $\ov{K}=\ker(s_{Y,\ov{Z}})$ under the natural homomorphism
$K\to \ov{K}$. Note that the exact sequence 
\begin{equation}\label{Z0-seq}
0\to C\to Z_0\to \ov{Z}\to 0
\end{equation}
gives rise to an exact sequence
$$0\to C\to K\to \ov{K}\to 0$$
Thus, we have an exact sequence of groups
$$1\to C\to G\to \ov{G}\to 1,$$
where $C$ is a central subgroup in $G$ and $\ov{G}$ is a finite Heisenberg group.
Thus, as in the proof of Theorem \ref{lag-inv-thm} we obtain an equivalence  
$$\FF(Y,\a)^{(Z_0,\b)}\cap\Coh(X)\simeq G-\rep_1.$$
Since $\ov{G}$ is a Heisenberg group, the desired equivalence $\Phi$ follows from 
Lemma \ref{center-lem}. Note that the natural functor 
$$\FF(Y,\a)^{(\ov{Z},\ov{\beta})}\to\FF(Y,\a)^{Z_0,\beta}$$
corresponds to the restriction under the homomorphism $G\to\ov{G}$. This implies
that $\Phi(S_{\ov{Z},\ov{\beta}})$ corresponds to the trivial one-dimensional representation of $C$.

Now considering the dual of the exact sequence 
\eqref{Z0-seq} we get (using the fact that $\ov{Z}$ is Lagrangian) the exact sequence
$$C^*\to X/\ov{Z}\to \hat{Z_0}\to\ldots$$
which gives an isomorphism of $Z'/\ov{Z}$ with $C^*$.
Thus, the pairing $\chi:Z'/\ov{Z}\times C^*\to\G_m$ is obtained as the canonical
pairing associated with the biextension $(\id\times j)^*\bar{\EE}$ of $X/\ov{Z}\times Z_0$,
where $\bar{\EE}$ is the biextension of $X/\ov{Z}\times \ov{Z}$ induced by $\EE$.
Explicitly, let 
$$t(C):\OO_{X/\ov{Z}\times C}\to  (\id\times j)^*\bar{\EE}|_{X/\ov{Z}\times C} \ \text{ and } 
t(Z'/\ov{Z}):\OO_{Z'/\ov{Z}\times Z_0}\to (\id\times j)^*\bar{\EE}|_{Z'/\ov{Z}\times Z_0}$$
be natural trivializations of restricted biextensions. Then
$$\chi\cdot t(C)|_{Z'/\ov{Z}\times C}=t(Z'/\ov{Z})|_{Z'/\ov{Z}\times C}.$$

The subgroup $Z'\sub X$ acts on $X$ by translations, and we have an isomorphism of
$1$-cocycles of $Y\times Z_0$ with values in $\PPic(X)$
$$L_{(y,z_0),x+z'}\simeq L_{(y,z_0),x}\ot\BB_{z',j(z_0)}^{-1}\ot\BB_{y,z'}^{-1}.
$$
Using the trivialization of $(\id\times j)^*\EE|_{Z'\times \ov{Z}}$ induced by $t(Z'/\ov{Z})$
we obtain an isomorphism
\begin{equation}\label{Y-Z0-Z'-cocycle-isom}
\tau^{-1}:L_{(y,z_0),x+z'}\rTo{\sim} L_{(y,z_0),x}\ot\BB_{j(z_0),z'}^{-1}\ot\BB_{y,z'}^{-1}\rTo{\sim}
L_{(y,z_0),x}\ot\BB_{y+j(z_0),z'}^{-1}.
\end{equation}
Thus, for fixed $z'\in Z'$
we are in the situation of Proposition \ref{cocycle-prop}(iii) with $H=Y\times Z_0$, $f:X\to X$ the
translation by $z'$, $M=\BB|_{X\times\{z'\}}$, the isomorphism $\tau_f=\tau$ given by
\eqref{Y-Z0-Z'-cocycle-isom} and the fixed point $x=0$. 
Note that the functor $T_{z'}$ on $\FF(Y,\a)$ sends $F$
to $t_{z'}^*F\ot M$, which is exactly the functor $\iota_{f,M}$ considered in Proposition
 \ref{cocycle-prop}(iii).
Let $G_{z'}$ be the central extension of $K$ by $\G_m$ with the underlying line bundle
$L|_{K\times\{z'\}}$. The diagram \eqref{fib-aut-diagram} 
gives in our case an isomorphism of $G$-representations
$$\fib_0(T_{z'}F)\simeq \a^*\fib_0(F)$$
for $F\in \FF(Y,\a)^{(Z_0,\b)}$, where the automorphism
$\a=\kappa^{-1}\circ\ov{\tau}:G\to G$ (identical on $\G_m\sub G$) is the composition
of isomorphisms of central extensions $\ov{\tau}:G\to G_{z'}$ and $\kappa^{-1}:G_{z'}\to G$
defined as follows. The isomorphism  
$\ov{\tau}:G\to G_{z'}$ is obtained by specializing \eqref{Y-Z0-Z'-cocycle-isom} to 
$(y,z_0)\in K$ and $x=0$, and
using the trivialization of $\BB_{0,z'}$.
On the other hand, choosing $(y',z'_0)$ such that $y'+j(z'_0)=z'$ we obtain another isomorphism
$\kappa^{-1}:G_{z'}\to G$ given by the isomorphism
\begin{align*}
&L_{(y,z_0),y'+j(z'_0)}=L_{(y,z_0),0}\ot\BB^{-1}_{y'+j(z'_0),j(z_0)}\ot\BB^{-1}_{y,y'+j(z'_0)}\simeq\\
&
L_{(y,z_0),0}\ot\BB^{-1}_{y',j(z_0)}\ot\BB^{-1}_{j(z_0),j(z'_0)}\ot\BB^{-1}_{y',y}\ot\BB^{-1}_{y,j(z'_0)}\simeq
\\
&L_{(y,z_0),0)}\ot\BB^{-1}_{y',y+j(z_0)}\ot\BB^{-1}_{y+j(z_0),j(z'_0)}\simeq
L_{(y,z_0),0)},
\end{align*}
where we used the symmetry of $\BB|_{Y\times Y}$ and of $\BB|_{\ov{Z}\times\ov{Z}}$ together
with trivializations of $\BB_{y',0}$ and of $\BB_{0,j(z'_0)}\simeq\BB_{j(z'_0),0}$. 
One can easily see from this that the restriction of $\a=\kappa^{-1}\circ\ov{\tau}$ to 
$C\times\G_m\sub G$ is given by $(c,\lambda)\mapsto (c,\chi_{z'}(c)\lambda)$,
where $\chi(z'):C\to \G_m$ is the character corresponding to $z'\mod\ov{Z}$.

\noindent
(ii) First, let us consider the case when $Y$ and $Z$ are transversal.
Set $\ov{Z}=j(Z_0)\sub X$. By Lemma \ref{conn-comp-lem}, $\ov{Z}$ is Lagrangian, so we can choose $\ov{\beta}$,
so that $(\ov{Z},\ov{\beta})$ is a Lagrangian pair. Now we have two completions of $Z_0$ to
an isotropic pair: $(Z_0,\beta|_{Z_0})$ and $(Z_0,j^*\ov{\beta})$. Hence,
$\beta|_{Z_0}\simeq j^*\ov{\beta}\ot \xi$ for some $\xi\in\hat{Z_0}$. Let us choose
$x\in X$ such that $\xi_{z_0}\simeq \EE_{j(z_0),x}$ and set
$$\beta'=\beta\ot \EE^{-1}|_{Z\times\{x\}}$$
Then $\beta'|_{Z_0}\simeq j^*\ov{\beta}$.
On the other hand, by Lemma \ref{beta-change-lem}, we have an equivalence
$$t(Z,x):\FF(Y,\a)^{(Z,\b)}\to\FF(Y,\a)^{(Z,\b')}$$
induced by $T_x$. This equivalence sends $S_{Z,\b}$ to
$T_x(S_{Z,\b})\simeq S_{Z,\b'}$. Thus, it is enough
to prove our statement with $\b$ replaced by $\b'$. In other words, we can
assume that $\b|_{Z_0}\simeq j^*\ov{\beta}$.

Let $G\to K=\ker(s_{Y,Z})$ (resp., $G_0\to K_0=\ker(s_{Y,Z_0})$)
be the central extensions by $\G_m$ appearing in the proof of Theorem \ref{lag-inv-thm} 
(resp., in the proof of (i)), so that we have equivalences 
\begin{equation}\label{Z-Z0-equiv-eq}
\FF(Y,\a)^{(Z,\b)}\cap\Coh(X)\simeq G-\rep_1 \ \text{ and }\ 
\FF(Y,\a)^{(Z_0,\b)}\cap\Coh(X)\simeq G_0-\rep_1.
\end{equation}
We have a commutative diagram of groups in which the horizontal
arrows are injective
\begin{diagram}
G_0&\rTo{}& G\\
\dTo{}&&\dTo{}\\
K_0&\rTo{}& K
\end{diagram}
Under the equivalences \eqref{Z-Z0-equiv-eq} the  natural functor
$$\FF(Y,\a)^{(Z,\b)}\cap\Coh(X)\to\FF(Y,\a)^{(Z_0,\b)}\cap\Coh(X)$$
corresponds to the restriction functor
$$G-\rep_1\to G_0-\rep_1\simeq C-\rep,$$
where $C=\ker(j|_{Z_0})$. 

We claim that the subgroup $Z'\sub X$
considered in (ii) coincides with $j(Z)$. Indeed, using the definition of $Z'$ and the fact
that $j:Z\to X$ is generalized Lagrangian we obtain a morphism of exact triangles
\begin{diagram}
Z & \rTo{}& X &\rTo{}& \Du(Z)&\rTo{}&Z[1]\\
\dTo{}&&\dTo{\id}&&\dTo{}&&\dTo{}\\
Z'&\rTo{}&X&\rTo{}&\hat{Z_0}&\rTo{}&Z'[1]
\end{diagram}
Since $H^0\Du(Z)\simeq \hat{Z_0}$, this implies the surjectivity of the map $Z\to Z'$,
induced by $j$, which proves our claim. 

Assume that $\Pi\simeq C^*$ is reduced.
Then the restriction of the Schr\"odinger representation $V_G$ of $G$ to $G_0$ corresponds under the equivalence $G_0-\rep_1\simeq C-\rep$ to a representation
$\bigoplus_{\chi\in C^*}\chi^{\oplus m_{\chi}}$.
Let 
$$V_G=\bigoplus_{\chi\in C^*}(V_G)_{\chi}$$
be the decomposition of $V_G$, viewed as a representation of $C$, into isotypic components.
It is well known that all $(V_G)_\chi$ have the same dimension (for example, this can be checked
by embedding $C$ into a maximal isotropic subgroup of $K$), hence
$$\dim (V_G)_\chi=\frac{\dim V_G}{|\Pi|}.$$
Using Lemma \ref{center-lem}, we see that 
$$m_\chi\dim V_{\ov{G}}= \dim (V_G)_\chi,$$
where $V_{\ov{G}}$ is the Schr\"odinger representation of $\ov{G}=G_0/C$. Thus,
we obtain the following formula for the multiplicities:
\begin{equation}\label{multiplicity-formula}
m_\chi=\frac{\dim V_G}{|\Pi|\cdot \dim V_{\ov{G}}}.
\end{equation}

Since $S_{\ov{Z},\ov{\b}}$ corresponds to the trivial representation of
$C$, the isomorphism \eqref{S-Z-Z0-coset-eq} will follow now from 
\eqref{translation-chi-isom}, once we show that
$m_\chi=N_Z$, i.e.,
$$\frac{\dim V_G}{\dim V_{\ov{G}}}=|\Pi|^{1/2}\cdot |\pi_0(Z)|^{1/2}$$
(note that $\Pi=j(Z)/j(Z_0)=\pi_0(j(Z))$). 
Equivalently, we have to check that
$$\frac{\deg(s_{Y,Z})}{\deg(s_{Y,j(Z_0)})}=|\Pi|\cdot |\pi_0(Z)|.$$
Consider the commutative diagram of isogenies
\begin{diagram}
Y\times Z_0 &\rTo{}& Y\times Z\\
\dTo{\id\times j}&&\dTo{s_{Y,Z}}\\
Y\times j(Z_0) &\rTo{s_{Y,j(Z_0)}}& X
\end{diagram}
where the top horizontal arrow is an embedding of the connected component of zero.
Note that 
\begin{equation}\label{C-Pi-eq}
\deg(j:Z_0\to j(Z_0))=|C|=|\Pi|
\end{equation}
due to duality between $\Pi$ and $C$.
Hence, this diagram leads to the following equality of degrees
$$\deg(s_{Y,Z_0})=\frac{\deg(s_{Y,Z})}{|\pi_0(Z)|}=|\Pi|\cdot \deg(s_{Y,j(Z_0)})$$
which gives the desired identity.

In the case when $\Pi$ is not reduced the representation of $C$ corresponding to $V_G|_{G_0}$
has a composition series, where the multiplicity $m_{\chi}$ of a character $\chi\in C^*(k)$ 
is given by a formula similar to \eqref{multiplicity-formula} but with $|\Pi|$ replaced by $\# \Pi(k)$.
This leads to the factor $\frac{|\Pi|}{\# \Pi(k)}$ in the multiplicities of the composition series
of $S_{Z,\b}$  in $\FF(Y,\a)^{(Z_0,\b)}\cap\Coh(X)$.

Now let us consider the case when $Y$ and $Z$ are not necessarily transversal.
As in the proof of Theorem \ref{lag-inv-thm}, using Lemma \ref{transv-lem} we choose
a Lagrangian pair $(Y',\a')$ transversal to $Z$, such that we have an equivalence
$\psi:\FF(Y,\a)_S\simeq\FF(Y',\a')_S$ compatible with $\bH$-action. 
Then $\psi$ induces an equivalence
$$\FF(Y,\a)^{(Z_0,\b)}\simeq \FF(Y',\a')^{(Z_0,\b)}$$
(and similar equivalences for $(Z,\b)$-invariants and for $(j(Z_0),\ov{\b})$-invariants).
As we have seen in the proof of Theorem \ref{lag-inv-thm},
$\psi(S_{j(Z_0),\ov{\b}}(Y,\a))$ is cohomologically
pure, so changing $\psi$ by $\psi[m]$ for appropriate $m\in\Z$
we can assume that $\psi$ induces an equivalence
$$\FF(Y,\a)^{(j(Z_0),\ov{\b})}\cap\Coh(X)\simeq \FF(Y',\a')^{(j(Z_0),\ov{\b})}\cap\Coh(X).$$
We claim that in this case $\psi$ also induces an equivalence
\begin{equation}\label{F-Y-Y'-Z0-equiv}
\FF(Y,\a)^{(Z_0,\b)}\cap\Coh(X)\simeq \FF(Y',\a')^{(Z_0,\b)}\cap\Coh(X).
\end{equation}
Indeed, we know that all simple objects of the finite length abelian category
$\FF(Y',\a')^{(Z_0,\b)}\cap\Coh(X)$ are of the form $T_xS$, where
$S=S_{j(Z_0),\ov{\b}}(Y',\a')$. But $\psi^{-1}(T_xS)=T_x\psi^{-1}(S)$ is a coherent
sheaf. Hence, $\psi^{-1}$ sends $\FF(Y',\a')^{(Z_0,\b)}\cap\Coh(X)$ to coherent sheaves.
Since the subcategory $\FF(Y',\a')^{(Z_0,\b)}\sub D^b(X)$ is compatible with the devissage with
respect to the standard $t$-structure, this easily implies that $\psi$
sends $\FF(Y,\a)^{(Z_0,\b)}\cap\Coh(X)$ to coherent sheaves and our claim follows.
Using the equivalence \eqref{F-Y-Y'-Z0-equiv}, we can transfer the 
decomposition of $S_{Z,\b}(Y',\a')$ in the category $\FF(Y',\a')^{(Z_0,\b)}\cap\Coh(X)$ to 
that of $S_{Z,\b}(Y,\a)$ in the category
$\FF(Y,\a)^{(Z_0,\b)}\cap\Coh(X)$.
\ed

\begin{cor}\label{gen-Lag-cor}
Let $(j:Z\to X,\b)$ be a generalized Lagrangian pair. Assume that $\cha(k)=0$.
Then 
$$S_{Z,\b}\simeq \bigoplus_{i=1}^N S_{j(Z_0),\b_i}$$
for some line bundles $\b_1,\ldots,\b_N$ on $j(Z_0)$ such that
$(j(Z_0),\b_i)$ are Lagrangian pairs.
\end{cor}

\Pf . This follows from Proposition \ref{lag-inv-prop}(ii) together with Lemma \ref{beta-change-lem}.
\ed

\begin{defi} Let $X=X_A$ and $Y=\{0\}\times\hat{A}$, so that we have $\FF(Y,\OO)\simeq D^b(A)$
(see Example \ref{main-rep-ex}). A coherent sheaf $F$ on $A$ is called a
{\it gLI-sheaf} if there exists a generalized Lagrangian pair $(Z,\b)$ such that
$F$ is $(Z,\b)$-invariant. We say that $F$ is an {\it LI-sheaf} if $Z$ can be chosen to
be a subvariety in $X_A$.
\end{defi}

We can show that for an LI-sheaf the subvariety $Z$ in the above definition
can be recovered as the stabilizer with respect to the action of $\bH$. Here is a slightly more
general result.

\begin{prop}\label{lag-max-prop} 
Let $(Y,\a)$ and $(Z,\b)$ be Lagrangian pairs for an ess-abelian variety $X$.
Assume that the condition $(\star)$ from Theorem \ref{lag-inv-thm} is satisfied 
and consider the generating object $S_{Z,\b}=S_{Z,\b}(Y,\a)\in\FF(Y,\a)^{(Z,\b)}\cap\Coh(X)$.
Then the subset of points $x\in X$ such that
$T_x(S_{Z,\b})\simeq S_{Z,\b}$ coincides with $Z$.
\end{prop}

\Pf . The same argument as in the proof of Theorem \ref{lag-inv-thm}
(based on Lemma \ref{transv-lem} and Theorem \ref{main-eq-thm}) shows that
it is enough to consider the case when $Y$ and $Z$ are transversal.
In this case $S_{Z,\b}$ is a vector bundle (see Theorem \ref{lag-inv-thm}).
Suppose $T_{x}(S_{Z,\b})\simeq S_{Z,\b}$ for some $x\in X$.
By definition of $(Z,\b)$-invariance, we have an isomorphism 
$$\b_z\ot T_z(p^*S_{Z,\b})\simeq p^*S_{Z,\b}$$
in $\FF(Y,\a)_Z\sub D^b(X\times Z)$, where $z=z^{un}\in Z(Z)$ is the universal point and
$p:X\times Z\to X$ is the projection. 
This leads to isomorphisms
$$\b_z\ot T_zT_{x}(p^*S_{Z,\b})\simeq \b_z\ot T_{x}T_z(p^*S_{Z,\b})\simeq p^*S_{Z,\b}$$
in $\FF(Y,\a)_Z$.
Since the commutator in the Heisenberg groupoid is given by the biextension $\EE$,
we deduce an isomorphism
$$\EE_{z,x}\ot p^*S_{Z,\b}\simeq p^*S_{Z,\b}$$
on $Z\times X$. Restricting to $Z\times \{0\}$ we get a trivialization of the line bundle 
$\EE_{z,x}$ on $Z$. Since $Z$ is Lagrangian, this implies that $x\in Z$.
\ed

\begin{cor}\label{uniqueness-cor} Keep the assumptions of Proposition \ref{lag-max-prop}.
Suppose that $S_{Z,\b}$ is invariant with respect to another Lagrangian pair
$(Z',\b')$. Then $Z'=Z$ and $\b'=\b$.
\end{cor}

\Pf. Proposition \ref{lag-max-prop} implies that $Z'\sub Z$, hence $Z'=Z$ (since $\dim Z=\dim Z'$).
By Lemma \ref{transv-lem}, we can assume that $Y$ is transversal
to $Z$, so that $S_{Z,\b}\in\FF(Y,\a)$ is a vector bundle (see Theorem \ref{lag-inv-thm}).
Now the isomorphism
$$\b_z\ot T_z(S_{Z,\b})\simeq \b'_z\ot T_z(S_{Z,\b})$$
on $Z\times X$ leads to 
$$\b_z^{-1}\ot\b'_z\ot S_{Z,\b}\simeq S_{Z,\b}.$$
Restricting to $Z\times\{0\}$ we deduce the triviality of $\b^{-1}\ot \b'$.
\ed

The formula for $S_{Z,\b}(Y,\a)$ from Theorem \ref{lag-inv-thm} has the following analog
in the non-transversal case.

\begin{prop}\label{non-transv-prop} 
Let $(Y,\a)$ be a Lagrangian pair, and let $(j:Z\to X,\b)$ be a generalized
Lagrangian pair. Assume that the condition $(\star)$ from Theorem \ref{lag-inv-thm} holds.
Assume in addition that the pull-backs of $\a$ and $\b$ to the
connected component of zero in $Y\times_X Z$ are isomorphic. Then
$S_{Z,\b}(Y,\a)$ is a direct summand in 
$$(s_{Y,Z})_*(\a^{-1}\boxtimes\b^{-1}\otimes (i\times j)^*\BB^{-1}),$$
where $i:Y\to X$ is the embedding map.
\end{prop}

\Pf . We can factor the homomorphism $s_{Y,Z}:Y\times Z\to X$ as a composition
$Y\times Z\rTo{q} X'\rTo{\iota} X$, where $q$ is surjective and $\iota$ is injective.
Now as in the proof of Theorem \ref{lag-inv-thm}, we obtain a structure of $1$-cocycle with respect
to the action of $Y\times Z$ on $X'$ on the line bundle $L$ on $(Y\times Z)\times X'$
given by the restriction of \eqref{L-cocycle}. Furthermore, we see that for an 
$L$-twisted $Y\times Z$-equivariant sheaf $F$ on $X'$
one has $\iota_*F\in \FF(Y,\a)^{(Z,\b)}$. The coherent sheaf
$$S=q_*(L|_{Y\times Z\times\{0\}})=q_*(\a^{-1}\boxtimes\b^{-1}\otimes (i\times j)^*\BB^{-1})$$
 on $X'$
has a natural $L$-twisted $Y\times Z$-equivariant structure coming from the structure of a $1$-cocycle
on $L$. Thus, by Theorem \ref{lag-inv-thm}, $S_{Z,\b}(Y,\a)$ is a direct summand of $S$
provided $S\neq 0$.

Now we observe that for $y+j(z)=0$, where $y\in Y$, $z\in Z$, we have an isomorphism
$$\BB_{y,j(z)}^{-1}\simeq \a_0^{-1}\ot\a_y\ot \a_{-y}.$$ 
Hence, we obtain the isomorphism
of the restriction of $L$ to $\ker(s_{Y,Z})\times\{0\}$ with
$\a_{-y}\ot \b^{-1}_z$. Under the isomorphism 
$Y\times_X Z\to \ker(s_{Y,Z}):(y,z)\mapsto (-y,z)$ this line bundle corresponds
to the difference between the restrictions of $\a$ and $\b$.
Thus, if $A$ is the connected component of zero in $\ker(s_{Y,Z})$, then our assumption on
$\a$ and $\b$ implies that $L|_{A\times\{0\}}\simeq\OO$, hence $S\neq 0$.
\ed

\section{Functors associated with Lagrangian correspondences}\label{corr-sec}

\subsection{Lagrangian correspondences for ess-abelian varieties}\label{lag-cor-sec}

\begin{defi} Let $(X,\BB_X)$ and $(Y,\BB_Y)$ be ess-abelian varieties. A {\it generalized Lagrangian correspondence} ({\it g-Lag-correspondence} for short) 
from $X$ to $Y$ is a generalized Lagrangian pair $(L,\a)$ for 
$(X\times Y, \BB_X^{-1}\boxtimes\BB_Y)$. In the case when
$L$ is a subvariety of $X\times Y$ (i.e., $(L,\a)$ is a Lagrangian pair), we say
that we have a {\it Lagrangian correspondence}.
Note that with every g-Lag-correspondence $(L,\a)$ from $X$ to $Y$ one can associate the 
{\it opposite g-Lag-correspondence} $(\si(L),\a^{-1})$ from $Y$ to $X$, where
$\si:X\times Y\to Y\times X$ is a natural isomorphism.
\end{defi}

\begin{ex}\label{Lag-corr-ex}
Let $f\in \Hom(X,Y)_{\Q}$ be an isomorphism between
$X$ and $Y$ in the category $\Ab_{\Q}$, so that we have an equality
\begin{equation}\label{Q-sym-isom}
\eta_X=\hat{f}\circ\eta_Y\circ f,
\end{equation}
where $\eta_X:X\to\hat{X}$ and $\eta_Y:Y\to\hat{Y}$ are the symplectic self-dualities
(in this situation we say that $f$ is {\it symplectic}).
Then similarly to Example \ref{Lag-ex}.1, we consider the graph of $f$
$$L(f)=\{(Nx, Nf(x))\ |\ x\in X\} \sub X\times Y,$$
where $N>0$ is an integer such that $Nf\in\Hom(X,Y)$. The equality
\eqref{Q-sym-isom} implies an isomorphism
$$(Nf\times Nf)^*\EE_Y\simeq \EE_X^{N^2}\simeq (N\id_X\times N\id_X)^*\EE_X$$
of biextensions of $X\times X$. Let $p_X:L(f)\to X$, $p_Y:L(f)\to Y$ and $\pi:X\to L(f)$
be the natural maps. Then the previous isomorphism can be rewritten as
$$(\pi\times\pi)^*(p_Y\times p_Y)^*\EE_Y\simeq(\pi\times\pi)^*(p_X\times p_X)^*\EE_X.$$
Hence, $(p_Y\times p_Y)^*\EE_Y\simeq (p_X\times p_X)^*\EE_X$, i.e., 
$L(f)$ is isotropic. Since $\dim L(f)=\dim X=\dim(X\times Y)/2$, by Lemma \ref{lag-lem}(ii), it
is Lagrangian. Conversely, it is easy to see that all
Lagrangian correspondences $L\sub X\times Y$, finite over $X$, are 
obtained by the above construction (in particular, they are automatically finite over $Y$).
\end{ex}

If $(L,\a)$ is a g-Lag-correspondence from
$X$ to $Y$ and $(M,\b)$ is a g-Lag-correspondence from $Y$ to $Z$
then we can try to define the composition $(M\circ L, \b\circ\a)$ by setting
$M\circ L=L\times_Y M$ and defining $\b\circ\a$ as the tensor product of the pull-backs
of $\a$ and $\b$ to $L\times_Y M$. Below we will give a 
sufficient condition for $(M\circ L,\b\circ\a)$
to be a g-Lag-correspondence from $X$ to $Z$ (see Corollary \ref{comp-cor}). 

It is convenient to consider the following more general setup.
Let $(X,\eta,\BB)$ be an ess-abelian variety.
Assume that $I\sub X$ is an isotropic abelian subvariety so that the composition 
$I\to X\stackrel{\eta}{\to} \hat{X}\to \hat{I}$ is zero. In other words, we have $I\sub I^{\perp}$,
where $I^{\perp}$ is the kernel of the composition $X\to\hat{X}\to\hat{I}$.
Then we can define the reduced ss-abelian variety $(\ov{X},\ov{\eta})$, where
$\ov{X}=I^{\perp}/I$ and $\ov{\eta}$ is induced by $\eta$.
Moreover, if we assume that the restriction $\BB|_{I\times I^{\perp}}$ is trivial then there is an induced
biextension $\ov{\BB}$ on $\ov{X}\times\ov{X}$, so that $(\ov{X},\ov{\BB})$ is an ess-abelian variety.

In this situation one can start with a generalized Lagrangian pair $(Y,\a)$ for $X$ and try to define the
corresponding Lagrangian pair for $\ov{X}$. 

\begin{prop}\label{red-prop} 
Let $(Y,\a)$ be a generalized Lagrangian pair for $X$ such that the natural map
$Y\to X/I^{\perp}$ is surjective. Set $\ov{Y}=Y\times_X I^{\perp}$ and let $\ov{\a}$ be the pull-back
of $\a$ to $\ov{Y}$. Then $(\ov{Y},\ov{\a})$ is a generalized Lagrangian pair for $\ov{X}$.
\end{prop}

\Pf . By assumption we have an exact sequence of commutative groups
$$0\to\ov{Y}\to Y\to X/I^{\perp}\to 0,$$
where the map $Y\to X/I^{\perp}$ is the composition of the natural maps $Y\to X$ and
$X\to X/I^{\perp}$. Dualizing and using the fact that $(Y,\a)$ is generalized
Largangian we obtain an
exact triangle in $D^b(\GG^{pr})$
$$I\to [Y\to X]\to \Du(\ov{Y})\to\ldots$$
where the first arrow is the composition $I\to X\to [Y\to X]$.
By the octahedron axiom we obtain also an exact triangle
$$Y\to X/I\to \Du(\ov{Y})\to\ldots.$$
In other words, $\Du(\ov{Y})$ is represented by the complex $[Y\to X/I]$.
It follows that the map $Y\to X/I$ has finite kernel. Now the cartesian square
\begin{diagram}
\ov{Y} &\rTo{}& Y\\
\dTo{}&&\dTo{}\\
\ov{X}&\rTo{}& X/I
\end{diagram}
leads to the exact triangle
$$\ov{Y}\to\ov{X}\to\Du(\ov{Y})\to\ldots$$
It is easy to deduce from this that the pair $(\ov{Y},\ov{\a})$ is generalized Lagrangian.
\ed

\begin{cor}\label{comp-cor} Let $(L,\a)$ (resp., $(M,\b)$)
be a g-Lag-correspondence from $X$ to $Y$ (resp., from $Y$ to $Z$).
Assume that the natural map $L\times M\to Y$ is surjective.
Then $(M\circ L,\b\circ\a)$ is a g-Lag-correspondence from $X$ to $Z$.
\end{cor}

\Pf . Consider the ess-abelian variety 
$$(X\times Y\times Y\times Z, \BB_X^{-1}\boxtimes\BB_Y\boxtimes\BB_Y^{-1}\boxtimes\BB_Z)$$
and the isotropic subvariety $I=0_X\times\De(Y)\times 0_Z$, where $\De:Y\to Y\times Y$ is the
diagonal embedding. To get the result we apply Proposition \ref{red-prop} to the generalized
Lagrangian pair $(L\times M,\a\boxtimes\b)$ for this ess-abelian variety. 
\ed

This result allows us to make the following

\begin{defi}\label{Lag-def} 
(i) For an ess-abelian variety $(X,\BB_X)$ we denote by
$\Lag(X)$ the set of isomorphism classes of g-Lag-correspondences $(L,\a)$ from $X$ to $X$
such that the projections $p_1,p_2:L\to X$ are surjective. The composition of correspondences
makes $\Lag(X)$ into a monoid, where the unit corresponds to $L=\De(X)\sub X\times X$ (and
trivial $\a$). 

\noindent
(ii)  For an ss-abelian variety $X$ 
we denote by $U(X,\Q)$ the group of symplectic automorphisms of $X$ in $\Ab_{\Q}$
(see Example \ref{Lag-corr-ex}).
\end{defi}

\begin{lem}\label{Lag-Sp-lem}
Associating with $(i:L\to X,\a)\in\Lag(X)$ the morphism in $\End(X)_{\Q}$
given by the correspondence
$i(L_0)\sub X\times X$ gives a surjective homomorphism of monoids
$$\pi:\Lag(X)\to U(X,\Q).$$
\end{lem}

\Pf . Since the projection $p_1\circ i:L_0\to X$ is an isogeny there exists a morphism
$\phi:X\to L_0$ such that 
$$i(\phi(x))=(nx,a(x))$$ 
for  some integer $n>0$ and some element $a\in\End(X)$.
By definition, $\pi(L)=a/n\in\End(X)_{\Q}$.
Now let $(j:M\to X,\b)$ be another element of $\Lag(X)$, and let $\psi:X\to M_0$ be a morphism
such that 
$$j(\psi(x))=(mx,b(x))$$
for some integer $m>0$ and some element $b\in\End(X)$, and so 
$\pi(M)=b/m$.
Then we have a morphism
$$X\to M\circ L=L\times_X M: x\mapsto (mnx, ma(x), ba(x))$$ 
which factors through the connected component of zero in $M\circ L$.
Thus,
$$\pi(M\circ L)=\frac{ba}{mn}=\pi(M)\circ\pi(L).$$
The homomorphism $\pi$ is surjective since the map $g\mapsto L(g)$ gives its set-theoretic
section (see Example \ref{Lag-corr-ex}).
\ed

Later we will need the following simple result about the composition of  
correspondences. 
For a g-Lag-correspondence $L\to X\times Y$ let us set
\begin{equation}\label{q-L-eq}
q(L)=\deg(L\to X)
\end{equation}
with the convention that this is zero if the projection $L\to X$ is not an isogeny.

\begin{lem}\label{q-lem} 
For g-Lag-correspondence $L\to X\times Y$ and
$M\to Y\times Z$ such that the map $L\times M\to Y$ is surjective one has
$$q(M\circ L)=q(L)\cdot q(M).$$
\end{lem}

\Pf . The projection $L\times_Y M\to X$ factors as the composition
$L\times_Y M\to L\to X$, so
$$q(M\circ L)=q(L)\cdot\deg(L\times_Y M\to L).$$
Now the cartesian square
\begin{diagram}
L\times_Y M &\rTo{} &M\\
\dTo{}&&\dTo{}\\
L&\rTo{}& Y
\end{diagram}
shows that $\deg(L\times_Y M\to L)=\deg(M\to Y)=q(M)$.
\ed

\subsection{LI-kernels and functors}\label{ker-fun-sec}

Recall that with every abelian variety $A$ we associate the standard ess-abelian variety 
$X_A=A\times\hat{A}$
equipped with a symplectic biextension $\EE_A=\BB_A\ot\si^*\BB_A^{-1}$,
where $\BB_A=p_{14}^*\PP$. The corresponding Heisenberg groupoid acts on
the category $D^b(A)$ (see Example \eqref{main-rep-ex}). 

\begin{defi}
(i) Let $A$ and $B$ be abelian varieties. Given a g-Lag-correspondence $(L,\a)$ 
from $X_A$ to $X_B$ we can consider the setup of Theorem \ref{lag-inv-thm}
for the Heisenberg groupoid associated with $(X_A\times X_B,\BB_A^{-1}\boxtimes\BB_B)$,
acting on the category $D^b(A\times B)$ that can be identified with $\FF(Y,\a)$ for 
$Y=\{0\}\times\hat{A}\times\{0\}\times\hat{B}\sub X_A\times X_B$ and $\a=\OO$
(see Remark \ref{F-Y-rem}.1). Therefore, we have the $t$-generator
$S_{L,\a}\in \Coh(A\times B)$ of $(L,\a)$-invariant objects in $D^b(A\times B)$
(see Theorem \ref{lag-inv-thm}).
We call $S_{L,\a}$ the {\it gLI-kernel defined by} $(L,\a)$ and denote
the corresponding {\it gLI-functor} by
$$\Phi_{L,\a}:=\Phi_{S_{L,\a}}:D^b(A)\to D^b(B).$$ 
In the case when $L$ is a subvariety of $X_A\times X_B$ we call $S_{L,\a}$ (resp.,
$\Phi_{L,\a}$) the {\it LI-kernel (resp., LI-functor) defined by} $(L,\a)$ 

\noindent
(ii) We say that a g-Lag-correspondence
$L\to X_A\times X_B=A\times\hat{A}\times B\times\hat{B}$ 
is {\it nondegenerate} if the projection $p_{AB}:L\to A\times B$ is surjective. 
\end{defi}

For a g-Lag-correspondence $L\to X_A\times X_B$ we will denote projections to
products of factors in $X_A\times X_B$ as
$p_A:L\to A$, $p_{A\hat{A}}:L\to A\times\hat{A}$, etc.

\begin{ex}\label{nondeg-ex} 
Let $L(g)$ be the Lagrangian correspondence from $X_A$ to itself 
associated with an element $g\in U(X_A,\Q)$ (see Example \ref{Lag-corr-ex}).
Let us write 
$$g=\left(\begin{matrix} a & b \\ c & d\end{matrix}\right), \text{ where } 
a\in\Hom(A,A)_{\Q}, \ b\in\Hom(\hat{A},A)_{\Q},\ 
c\in\Hom(A,\hat{A})_{\Q},\ d\in\Hom(\hat{A},\hat{A})_{\Q}.$$ 
Then the correspondence $L(g)$ is nondegenerate if and only if
the map in $\Ab_\Q$
$$A\times\hat{A}\to A\times A: (x,\xi)\mapsto (x,ax+b\xi)$$ 
is an isomorphism. Equivalently, $b$ should be invertible in $\Ab_\Q$. 
\end{ex}

The proof of the following lemma is straightforward and is left to the reader.

\begin{lem}\label{T-A-B-lem}
Under the natural identification of $D^b(A\times B)$ with
$\FF(Y,\OO)$, where $Y=\{0\}\times\hat{A}\times\{0\}\times\hat{B}\sub X_A\times X_B$
one has for $K\in D^b(A\times B)$
$$T_{(a,\xi,b,\eta)}(K)=t_{(a,b)}^*K\ot (\PP_{-\xi}\boxtimes\PP_{\eta})=
(\bbT^A_{a,\xi})^{-1}\circ_A K\circ_B\bbT^B_{b,\eta},$$
where $a\in A$, $\xi\in\hat{A}$, $b\in B$, $\eta\in\hat{B}$. Here $\bbT^A_x$ are the kernels
\eqref{bbT-eq} and $(\bbT^A_{a,\xi})^{-1}=\PP_{a,\xi}\ot \bbT^A_{-a,-\xi}$.
\end{lem}

This lemma allows to rewrite the condition of $(L,\a)$-invariance for $K\in D^b(A\times B)$
as an isomorphism
$$K\simeq \a_l\ot \left((\bbT^A_{p_{A\hat{A}}(l)})^{-1}\circ_A K\circ_B
\bbT^B_{p_{B\hat{B}}(l)}\right),$$
for $l\in L$, or equivalently,
\begin{equation}\label{intertw-eq}
\bbT^A_{p_{A\hat{A}}(l)}\circ_A K\simeq \a_l\ot K\circ_B\bbT^B_{p_{B\hat{B}}(l)}.
\end{equation}
Hence, the gLI-functor associated with $(L,\a)$ satisfies the following "intertwining"
isomorphisms involving Heisenberg groupoids actions on $D^b(A)$ and $D^b(B)$:
$$\Phi_{L,\a}\circ T_{p_{A\hat{A}}(l)}\simeq \a_l\otimes T_{p_{B\hat{B}}(l)}\circ\Phi_{L,\a}$$
for $l\in L$. 

In the case when both projections $L\to X_A$ and $L\to X_B$ are surjective, i.e.,
$(L,\a)$ is an element of $\Lag(X_A)$ (see Definition \ref{Lag-def}(i)), we can use
\eqref{intertw-eq} to move the kernels of the form $\bbT^A_x$ through $K$. We
record this observation for future use in the next lemma.

\begin{lem}\label{K-T-lem} 
Let $(L,\a)$ be a g-Lag-correspondence in $\Lag(X_A)$ and let $K\in D^b(A\times A)$
be an $(L,\a)$-invariant kernel. For any $x\in X_A$ there exists $x'\in X_A$ and
$x''\in X_A$ such that
$$K\circ_A \bbT^A_x=\bbT^A_{x'}\circ_A K \text{ and} \ \bbT^A_x\circ_A K=K\circ_A \bbT^A_{x''}$$
in $D^b(A\times A)$.
\end{lem}

From Theorem \ref{lag-inv-thm} and Proposition \ref{non-transv-prop} we get the following
(almost) explicit formulas for $S_{L,\a}$.

\begin{lem} (i) If $(L,\a)$ is a nondegenerate g-Lag-correspondence from $X_A$ to $X_B$
then $S_{L,\a}$ is a vector bundle on $A\times B$ given by
$$S_{L,\a}\simeq p_{AB*}\bigl(\a^{-1}\ot p_{A\hat{A}}^*\PP^{-1}\ot p_{B\hat{B}}^*\PP\bigr)^I,$$
where $p_{A\hat{A}}$ and $p_{B\hat{B}}$ are projections from $L$ to $X_A=A\times\hat{A}$ and
$X_B=B\times\hat{B}$, respectively; $I$ is a Lagrangian subgroup in $G:=\ker(p_{AB})$ with respect
to the central extension of $G$ associated with $\a|_G$.

\noindent
(ii) Let $(L,\a)$ be an arbitrary g-Lag-correspondence from $X_A$ to $X_B$.
Assume that the restriction of $\a$ to the connected component of zero in $\ker(p_{AB})$
is trivial. Then
$S_{L,\a}\in \Coh(A\times B)$ is a direct summand in
$$p_{AB*}\bigl(\a^{-1}\ot p_{A\hat{A}}^*\PP^{-1}\ot p_{B\hat{B}}^*\PP\bigr).$$
\end{lem}

\Pf . (i) Applying Theorem \ref{lag-inv-thm} to the generalized Lagrangians
$Z=L\to X_A\times X_B$ and $Y=\{0\}\times\hat{A}\times\{0\}\times\hat{B}\sub X_A\times X_B$, we
obtain that $S_{L,\a}$ as an object of $\FF(Y,\OO)$ is given by
$$S_{L,\a}\simeq q_*\bigl(\a(l)^{-1}\ot\PP_{p_A(l),\xi}\ot\PP^{-1}_{p_B(l),\eta}\bigr)^I,$$
where $(l,\xi,\eta)\in L\times \hat{A}\times\hat{B}$, and $q$ is the map 
$$L\times\hat{A}\times\hat{B}\to A\times\hat{A}\times B\times\hat{B}:
(l,\xi,\eta)\mapsto (p_A(l),p_{\hat{A}}(l)+\xi,p_B(l),p_{\hat{B}}(l)+\eta).$$
The identification of $\FF(Y,\OO)$ with $D^b(A\times B)$ is given by the restriction to
$A\times\{0\}\times B\times\{0\}\sub X_A\times X_B$ (see Example \ref{main-rep-ex}). 
Since $p_{AB}:L\to A\times B$ is surjective, the map $q$ is also surjective.
Therefore, we can use the base change
formula to get the required expression for $S_{L,\a}$ as an object of $D^b(A\times B)$.

\noindent (ii) This follows from Proposition \ref{non-transv-prop} by a similar argument.
Note that we can still use the base change formula since the image of $q$ is transversal
to the subvariety $A\times\{0\}\times B\times\{0\}\sub X_A\times X_B$.
\ed

\begin{exs}\label{LI-ex} 
1. According to \cite{O} any equivalence between $D^b(A)$ and $D^b(B)$ appears
as an LI-functor associated with the graph $L(f)$ of a symplectic isomorphism
$f:X_A\simeq X_B$ (see Example \ref{Lag-corr-ex}). More precisely, to construct such an
equivalence one uses the natural equivalence of $D^b(A)$ with
the representation $\FF(Y,\a)$ of the Heisenberg groupoid $\bH(X_B)$ of $X_B$ associated with the
Lagrangian pair $(Y,\a)$ that corresponds to $\{0\}\times\hat{A}$ under the isomorphism $f$.
Now Theorem \cite[Thm.\ 4.3]{P-sym} gives an equivalence of $\FF(Y,\a)$ with $D^b(B)$.
Note that the above construction depends on a choice of an extension of $f$ to an equivalence
of Heisenberg groupoids $\bH(X_A)\simeq\bH(X_B)$, and
the obtained equivalence $D^b(A)\simeq D^b(B)$ is compatible with the action
of these groupoids. This implies that the corresponding
kernel $K$ on $A\times B$ belongs to $D^b(A\times B)^{(L(f),\b)}$, where $(L(f),\b)$ is some Lagrangian
correspondence extending the graph $L(f)$. Furthermore, by \cite[Prop.\ 3.2]{O}, $K$ is cohomologically
pure. Hence, by Theorem \ref{lag-inv-thm}, $K$ is a direct sum of several copies of
$S_{L(f),\b}$, which implies that $K\simeq S_{L(f),\b}$ (since the corresponding functor is an
equivalence).
For example, the Fourier-Mukai transform $\SS:D^b(A)\to D^b(\hat{A})$ is
associated with the Lagrangian 
$\Ga_{\SS}\sub X_A\times X_{\hat{A}}=A\times\hat{A}\times\hat{A}\times A$
consisting of $(x,\xi,\xi,-x)$, where $x\in A$, $\xi\in\hat{A}$.
Note that in the case $A=B$ autoequivalences corresponding to symplectic automorphisms of 
$X_A$ were also considered in \cite{P-thesis} and \cite{Mukai-spin}.

\noindent 2. 
Let $f:A\to B$ be a homomorphism. Then with $f$ we can associate a Lagrangian correspondence
from $X_A$ to $X_B$ by setting $L=A\times\hat{B}$, where the map $L\to B$ is induced by $f$
and the map $L\to \hat{A}$ is induced by $\hat{f}$. In this case we can take $\a$ to be trivial.
Then $L$-invariants in $D^b(A\times B)$ are generated by the structure sheaf of the graph of $f$.
The corresponding LI-functor $D^b(A)\to D^b(B)$ is the derived push-forward
$Rf_*$. The functor corresponding to the opposite Lagrangian is the pull-back
$Lf^*:D^b(B)\to D^b(A)$.
\end{exs}

%

The following proposition shows that the adjoint functors to gLI-functors are also
gLI-functors.

\begin{prop}\label{dual-prop}
Let $(L,\a)$ be a g-Lag-correspondence from $X_A$ to $X_B$, and
consider the permutation maps
$$\si:A\times B\to B\times A \text{ and } \ \si_X:X_A\times X_B\to X_B\times X_A$$
Then
the objects $S_{L,\a}$ and $\si^*S_{\si_X(L),\a^{-1}}$ in $D^b(A\times B)$ are dual up to a shift.
Hence, the functors $\Phi_{L,\a}:D^b(A)\to D^b(B)$ and $\Phi_{\si L,\a^{-1}}:D^b(B)\to D^b(A)$ 
are adjoint up to a shift.
\end{prop}

\Pf . Consider the duality functor 
$$D: D^b(A\times B)\to D^b(A\times B): K\mapsto R\und{\Hom}(K,\OO).$$
We have
$$T_{(a,\xi,b,\eta)}(D(K))\simeq D(T_{(a,-\xi,b,-\eta)}(K))$$
for $(a,\xi,b,\eta)\in X_A\times X_B$. On the other hand, for $K'\in D^b(B\times A)$ we have
$$T_{(a,-\xi,b,-\eta)}(\si^*K')\simeq \si^*(T_{(b,\eta,a,\xi)}K').$$
Combining these isomorphisms we obtain
$$T_u(D(\si^*K'))\simeq D(\si^*(T_{\si_X(u)}K')),$$
where $u\in X_A\times X_B$. This shows that $D(\si^*S_{\si_X(L),\a^{-1}})$ has
an $(L,\a)$-invariant structure.
\ed

Our main goal is to calculate convolutions of gLI-kernels (and hence
compositions of gLI-functors). We will need the following technical result for this. 

\begin{lem}\label{coh-line-bun-lem} 
Let $X$ be a commutative group scheme equipped with homomorphisms $f:X\to A$, $g:X\to\hat{A}$,
where $A$ is an abelian variety. Assume that both $f$ and $g$ are surjective with finite kernel.
Assume also that $\LL$ is a line bundle on $X$ such that
\begin{equation}\label{biext-f-g-eq}
\La(\LL)\simeq (f\times g)^*\PP_A.
\end{equation}
Then the restriction of $\LL$ to every connected component of $X$ is nondegenerate
of the same index $i_{\LL}$, so $H^i(X,\LL)=0$ for $i\neq i_{\LL}$. Furthermore,
$$(\dim H^{i_{\LL}}(X,\LL))^2=|\ker(f)|\cdot|\ker(g)|  \text{ and}$$
$$i_{\LL}=i(\hat{f_0}\circ g_0), \text{ where } f_0=f|_{X_0}, g_0=g|_{X_0}$$
(here we use the notation of Sec.\ \ref{index-sec}).
\end{lem}

\Pf . Set $G=\ker(g)$, $H=\ker(f)$, and let $G_0=G\cap X_0$, $H_0=H\cap X_0$. Then we have
exact sequences
$$0\to G_0\to X_0\rTo{g_0} \hat{A}\to 0$$
$$0\to H_0\to X_0\rTo{f_0} A\to 0.$$
Also, we have $\pi_0(X)\simeq G/G_0\simeq H/H_0$.
The isomorphism \eqref{biext-f-g-eq} implies that for every point $x\in H$ one has
$t_x^*\LL\simeq\LL$. Since $H$ surjects onto $\pi_0(X)$, this implies that the restrictions of
$\LL$ to all connected components of $X$ are obtained from $\LL|_{X_0}$ by translations.
Hence, 
\begin{equation}\label{H-X-0-eq}
H^*(X,\LL)\simeq H^*(X_0,\LL|_{X_0})^{\oplus |\pi_0(X)|}.
\end{equation}
On the other hand, the isomorphism \eqref{biext-f-g-eq} induces an isomorphism
$$\La(\LL|_{X_0})\simeq (f_0\times g_0)^*\PP_A.$$
Therefore, the homomorphism $\phi:X_0\to\widehat{X_0}$ associated with $\LL|_{X_0}$
is equal to the composition 
$$X_0\rTo{f_0} A\rTo{\widehat{g_0}} \widehat{X_0}.$$
Therefore, the group $K:=K(\LL|_{X_0})=\ker(\phi)$ fits into an exact sequence
$$0\to H_0\to K\to G_0\to 0.$$
Hence, $K$ is finite, i.e., the line bundle $\LL|_{X_0}$ is nondegenerate of some index $i=i_{\LL}$.
It remains to observe that
$$|K|=|H_0|\cdot |G_0|=\frac{|H|\cdot |G|}{|\pi_0(X)|^2},$$
and so from \eqref{H-X-0-eq} we obtain
$$\dim H^i(X,\LL)=|\pi_0(X)|\cdot |K|^{1/2}=|H|^{1/2}\cdot |G|^{1/2}.$$
The formula for the index $i_{\LL}$ follows from the fact that
$\phi_{\LL}=\hat{f_0}\circ g_0$.
\ed

In the next proposition we calculate convolutions of gLI-kernels under some technical nondegeneracy
assumptions.

\begin{prop}\label{convolution-prop} 
(i) Let $(L,\a)$ (resp., $(M,\b)$)
be a g-Lag-correspondence from $X_A$ to $X_B$ (resp., from $X_B$ to $X_C$).
Assume that the natural homomorphism $L\times M\to X_B$ (restricting to the given ones
on $L$ and $M$) is surjective. 
Then for an $(L,\a)$-invariant object $F\in D^b(A\times B)$ and an $(M,\b)$-invariant
object $G\in D^b(B\times C)$ the convolution $F\circ_B G$ (see Sec.\ \ref{conv-sec})
can be equipped with a
$(M\circ L,\b\circ\a)$-invariance structure.

\noindent
(ii) Assume in addition that the g-Lag-correspondences $L$, $M$ and $M\circ L$ are nondegenerate.
Then one has an isomorphism in $D^b(A\times C)$
$$S_{L,\a}\circ_B S_{M,\b}\simeq S_{M\circ L,\b\circ\a}[\la(L,M)],$$
and hence an isomorphism of functors
$$\Phi_{M,\b}\circ\Phi_{L,\a}\simeq \Phi_{M\circ L,\b\circ\a}[\la(L,M)].$$
Here the integer $\la(L,M)\le 0$ is defined as follows. 
Let us set 
$$F=\ker(L\times_B M\to A\times C:(l,m)\mapsto (p_A(l),p_C(m))).$$
Consider the maps 
$$f: F\to B \text{ and} \ \ g: F\to \hat{B},$$
where $f$ is the restriction by the projection to $B$ and $g$ is the restriction of the map
$$L\times_B M\to \hat{B}:(l,m)\mapsto p_{\hat{B}}(l)-p_{\hat{B}}(m).$$
Set $f_0=f|_{F_0}$ and $g_0=g|_{F_0}$.
Then both $f_0$ and $g_0$ are isogenies, $\hat{f_0}\circ g_0$ is symmetric and we set
\begin{equation}\label{la-L-M-eq}
\la(L,M)=-i(\hat{f_0}\circ g_0).
\end{equation}
\end{prop}

\Pf . (i) For $\FF\in D^b(A\times B)$, $\GG\in D^b(B\times C)$ let us set
$$\FF\diamond \GG=p_{12}^*\FF\ot^{{\mathbb L}} p_{23}^*\GG\in F^b(A\times B\times C).$$
We have
$$(t^*_{(a,b)}\FF)\diamond(t^*_{(b,c)}\GG)\simeq t^*_{(a,b,c)}(\FF\diamond \GG) \text{ and}$$
$$(\FF\ot (\PP_{-\xi}\boxtimes\PP_{\eta}))\diamond
(\GG\ot (\PP_{-\eta}\boxtimes\PP_{\rho}))=(\FF\diamond \GG)\ot (\PP_{-\xi}\boxtimes\OO_B\boxtimes\PP_\rho),
$$
where $(a,b,c)\in A\times B\times C$ and $(\xi,\eta,\rho)\in \hat{A}\times\hat{B}\times\hat{C}$.
Hence, for $(l,m)\in L\times_{X_B}M$, $\FF\in D^b(A\times B)^{(L,\a)}$ and 
$\GG\in D^b(B\times C)^{(M,\b)}$
we obtain an isomorphism
$$\FF\diamond \GG\simeq T_l(\FF)\diamond T_m(\GG)\simeq
(\bbT^A_{p_{A\hat{A}}(l)})^{-1}\circ_A t_{b(l,m)}^*(\FF\diamond \GG)\circ_C\bbT^C_{p_{C\hat{C}}(m)},$$
where $b(l,m)=p_B(l)=p_B(m)$. Consider
the extension of the canonical map $L\times_{X_B} M\to X_A\times X_C$
to a map 
$i:L\times_{X_B} M\to X_A\times B\times X_C\sub X_A\times X_B\times X_C$ with the $B$-component
given by $b(l,m)$. We can view in this way $(M\circ L,\b\circ\a)$ 
as a generalized isotropic pair for $X_A\times X_B\times X_C$
equipped with the (degenerate) biextension $\BB_A^{-1}\boxtimes\OO_{X_B}\boxtimes\BB_C$. 
The above calculation shows that $\FF\diamond \GG$ has an 
$(M\circ L,\b\circ\a)$-invariance structure. This immediately leads to
the required invariance structure on the push-forward of $\FF\diamond \GG$ to $D(A\times C)$,
i.e., on $\FF\circ_B \GG$.

\noindent 
(ii) Let $p_{AB}:L\to A\times B$, $p_{BC}:M\to B\times C$ and $p_{AC}:L\circ M\to A\times C$
denote the natural projections.
Recall that 
$$\wt{S}_{L,\a}:=p_{AB*}(\a^{-1}\ot p^*_{A\hat{A}}\PP^{-1}\ot p^*_{B\hat{B}}\PP)\simeq
V_L\ot S_{L,\a},$$
$$\wt{S}_{M,\b}:=p_{BC*}(\b^{-1}\ot p^*_{B\hat{B}}\PP^{-1}\ot p^*_{C\hat{C}}\PP)\simeq
V_M\ot S_{M,\b},$$
where $V_L$ and $V_M$ are vector spaces of ranks $|\ker(p_{AB})|^{1/2}$ and
$|\ker(p_{BC})|^{1/2}$, respectively. Also, $S_{M\circ L,\b\circ\a}$ is a vector bundle
of rank $|\ker(p_{AC})|^{1/2}$.
Thus, by part (i) and Theorem \ref{lag-inv-thm}, it suffices to check that 
$$\wt{S}_{L,\a}\circ_B \wt{S}_{M,\b}\simeq \VV[-i(\hat{f_0}\circ g_0]$$ 
where $\VV$ is a vector bundle on $A\times C$ of rank 
$$|\ker(p_{AC})|^{1/2}\cdot |\ker(p_{AB})|^{1/2}\cdot |\ker(p_{BC})|^{1/2}.$$
Using the commutative diagram with
cartesian squares 
\begin{equation}\label{big-diagram-eq}
\begin{diagram}
&&&&L\times_B M\\
&&&\ldTo{}&&\rdTo{}\\
&&L\times C&&&&A\times M\\
&\ldTo{}&&\rdTo{}&&\ldTo{}&&\rdTo{}\\
L&&&&A\times B\times C&&&&M\\
&\rdTo{p_{AB}}&&\ldTo{}&&\rdTo{}&&\ldTo{p_{BC}}\\
&&A\times B&&&&B\times C\\
\end{diagram}
\end{equation}
and applying base change and projection formulas
we find an isomorphism
\begin{equation}\label{composition-sheaf-eq}
\begin{array}{l}
\wt{S}_{L,\a}\circ_B \wt{S}_{M,\b}\simeq\\
p_{AB*}\bigl(\a^{-1}\ot p_{A\hat{A}}^*\PP^{-1}\ot p_{B\hat{B}}^*\PP\bigr)\circ_B
p_{BC*}\bigl(\b^{-1}\ot p_{B\hat{B}}^*\PP^{-1}\ot p_{C\hat{C}}^*\PP\bigr)\simeq
p'_{AC*}(\LL),
\end{array}
\end{equation} 
where $p'_{AC}:L\times_B M\to A\times C$ is the projection and
$$\LL=\a^{-1}\ot\b^{-1}\ot p_{A\hat{A}}^*\PP^{-1}\ot p_{C\hat{C}}^*\PP^{-1}
\ot p_{B\hat{B}\hat{B}}^*(\id_B\times\de)^*\PP$$
with $\de:\hat{B}\times\hat{B}\to\hat{B}$ given by $\de(\xi_1,\xi_2)=\xi_1-\xi_2$.
It is easy to see that
$$\La(\LL)_{(l_1,m_1),(l_2,m_2)}\simeq \PP^{-1}_{p_A(l_2),p_{\hat{A}}(l_1)}\ot
 \PP_{p_C(m_2),p_{\hat{C}}(m_1)}\ot \PP_{p_B(l_2),p_{\hat{B}}(l_1)-p_{\hat{B}}(m_1)}.$$
Let us consider the fiber $F=\ker(p'_{AC}:L\times_B M\to A\times C)$.
Then the above formula specializes to
$$\La(\LL)|_{F\times F}\simeq (f\times g)^*\PP.$$
The diagram \eqref{big-diagram-eq} implies that the projection $L\times_B M\to A\times B\times C$
is surjective with a finite kernel $G$, where 
$$|G|=\deg(p_{AB})\cdot \deg(p_{BC})=|\ker(p_{AB})|\cdot |\ker(p_{BC})|.$$ 
Hence, the projection $f:F\to B$ is also surjective with the kernel isomorphic to $G$. 
On the other hand, the subgroup
$H=\ker(g)\sub F$ can be identified with 
$$(L\times_{X_B} M)\cap F=\ker(p_{AC}:M\circ L\to A\times C),$$
where we view $M\circ L=L\times_{X_B} M$ as a subgroup in $L\times_B M$.
Recall that by assumption, $p_{AC}$ is surjective with finite kernel, so $g:F\to \hat{B}$
is surjective with finite kernel.
Now Lemma \ref{coh-line-bun-lem} implies that
$p'_{AC*}(\LL)$ is of the form $\VV[-i(\hat{f_0}\circ g_0]$ for some vector bundle $\VV$ of rank 
$$|G|^{1/2}\cdot |H|^{1/2}=|\ker(p_{AC})|^{1/2}\cdot |\ker(p_{AB})|^{1/2}\cdot |\ker(p_{BC})|^{1/2},$$
as required.
\ed


To get rid of the nondegeneracy assumption in Proposition \ref{convolution-prop}
we will use twisting by autoequivalences. Namely, we will use the fact that in the
case when $L=L(f)\sub X_A\times X_A$ is a graph of a symplectic automorphism $f:X_A\to X_A$
the corresponding functor $\Phi_{L(f),\a}$ is an equivalence (see Example \ref{LI-ex}.1).

\begin{lem}\label{auteq-conv-lem} 
Let $(L,\a)$ be a g-Lag-correspondence from $X_A$ to $X_B$, and let
$f:X_A\to X_A$ (resp., $g:X_B\to X_B$) be a symplectic automorphism. Let also
$(L(f),\a_f)$ (resp., $(L(g),\a_g)$) be a Lagrangian correspondence from $X_A$ to $X_A$
(resp. from $X_B$ to $X_B$) extending the graph of $f$ (resp., $g$).
Then
$$S_{L(f),\a_f}\circ_A S_{L,\a}\circ_B S_{L(g),\a_g}\simeq S_{L(g)\circ L\circ L(f),\a_g\circ\a\circ\a_f}[i]$$
for some $i\in \Z$, and
\begin{equation}\label{g-L-f-eq}
L(g)\circ L\circ L(f)=(f^{-1}\times g)(L).
\end{equation}
\end{lem}

\Pf . Note that \eqref{g-L-f-eq} follows immediately from the definition of the composition of correspondences.
Let us set for brevity $S_f=S_{L(f),\a_f}$, $S_g=S_{L(g),\a_g}$, $S=S_{L,\a}$ and
$S'=S_{L(g)\circ L\circ L(f),\a_g\circ\a\circ\a_f}$. Let $S_{f^{-1}}$ denote the kernel
of the inverse autoequivalence to the one defined by $S_f$, so that
$$S_{f^{-1}}\circ_A S_f\simeq\De_*\OO_A.$$
Then $S_{f^{-1}}$ has a $(L(f^{-1}),\si^*\a_f^{-1})$-invariance structure, where $\si:L(f^{-1})\to L(f)$
is the isomorphism induced by the permutation of factors in $X_A\times X_A$,
so $S_{f^{-1}}$ is isomorphic to $S_{L(f^{-1}),\si^*\a_f^{-1}}[m]$ for some $m\in\Z$.
Similarly, for $S_{g^{-1}}=S_{L(g^{-1}),\si^*\a_g^{-1}}[n]$, where $n\in\Z$, we have
$$S_g\circ_B S_{g^{-1}}\simeq\De_*\OO_B.$$
Now by Proposition \ref{convolution-prop}(i), the object
$$P:=S_f\circ_A S\circ_B S_g\in D^b(A\times B)$$
has a structure of
an object of $D^b(A\times A)^{(L(g)\circ L\circ L(f),\a_g\circ\a\circ\a_f)}$.
Similarly the object
$$Q:=S_{f^{-1}}\circ_A S'\circ_B S_{g^{-1}}$$
has a structure of an object of $D^b(A\times B)^{(L,\a)}$.
Let $[a,b]$ (resp., $[s,t]$) be the cohomological amplitude of $P$ 
(resp., $Q$). Then by Theorem \ref{lag-inv-thm}, we have
an exact triangle
$$\tau_{<b}P\to P\to V\ot S'[-b]\to\ldots,$$
where $V$ is a vector space and $\tau_{<b}P$ is a successive extension of
$S'[-i]$ with $i<b$. Applying the convolution with $S_{f^{-1}}$ on the left and with
$S_{g^{-1}}$ on the right we obtain an exact triangle
$$S_{f^{-1}}\circ _A\tau_{<b}P\circ_B S_{g^{-1}}\to S\to V\ot Q[-b]\to\ldots,$$
where $S_{f^{-1}}\circ _A\tau_{<b}P\circ_B S_{g^{-1}}$ is a successive extension of
$Q[-i]$ with $i<b$.
It follows that 
$$\und{H}^{b+t}(S)\simeq V\ot \und{H}^t(Q)\neq 0.$$
A similar argument shows that $\und{H}^{a+s}(S')\neq 0$.
Hence, $b+t=a+s=0$ which implies that $b=a$ and $s=t=-a$.
Thus, we have
$$P\simeq V\ot S'[-a] \ \text{ and } \ Q\simeq W\ot S[a]$$
for some vector spaces $V$ and $W$. Therefore, we obtain
$$S\simeq S_{f^{-1}}\circ_A P\circ_B S_{g^{-1}}\simeq V\ot Q[-a]\simeq V\ot W\ot S,$$
which implies that $V$ and $W$ are one-dimensional. 
\ed

Now we are ready to prove our main result about the convolution of gLI-kernels.

\begin{thm}\label{main-conv-thm} 
Let $(L,\a)$ (resp., $(M,\b)$)
be a g-Lag-correspondence from $X_A$ to $X_B$ (resp., from $X_B$ to $X_C$).
Assume that the natural homomorphism $L\times M\to X_B$ is surjective. 
Then one has an isomorphism in $D^b(A\times C)$
$$S_{L,\a}\circ_B S_{M,\b}\simeq S_{M\circ L,\b\circ\a}[i] \ \text{ for some } i\in\Z,$$
and hence an isomorphism of functors
$$\Phi_{M,\b}\circ\Phi_{L,\a}\simeq \Phi_{M\circ L,\b\circ\a}[i].$$
\end{thm}

\Pf . The idea is to reduce to the case when the correspondences $L$, $M$ and $M\circ L$
are nondegenerate, considered in Proposition \ref{convolution-prop}(ii). Let
$\phi_A:\hat{A}\to A$, $\phi_B:\hat{B}\to B$ and $\phi_C:\hat{C}\to C$ be symmetric isogenies
associated with some ample line bundles on $\hat{A}$, $\hat{B}$ and $\hat{C}$.
For any $n\in\Z$ let us define the symplectic automorphism $f_A(n)$ of $X_A=A\times\hat{A}$ by
$$f_A(n)=\left(\begin{matrix} \id_A & n\phi_A \\ 0 & \id_{\hat{A}} \end{matrix}\right).$$
Similarly, we define symplectic automorphisms $f_B(n)$ of $X_B$ and $f_C(n)$ of $X_C$.
We claim that for all $n$ except for a finite number of values the correspondence
$(f_A(n)\times f_B(n))(L)$ is nondegenerate. Indeed, by definition, nondegeneracy means transversality
to the Lagrangian subvariety $\{0\}\times \hat{A}\times \{0\}\times\hat{B}\sub X_A\times X_B$.
Thus, we need $n$ such that $L$ is transversal to 
$$(f_A(-n)\times f_B(-n))(\{0\}\times \hat{A}\times \{0\}\times\hat{B})=
\Ga(-n(\phi_A\times\phi_B)),$$
and our claim follows from Lemma \ref{transv-lem}(ii). Similarly, for generic $n$ the correspondences
$(f_B(n)\times f_C(n))(M)$ and $(f_A(n)\times f_C(n))(M\circ L)$ are nondegenerate.
Let us choose such $n$ and extend $f_A=f_A(n)$, $f_B=f_B(n)$ and $f_C=f_C(n)$ to some Lagrangian
correspondences $(L(f_A),\a_1)$, $(L(f_B),\a_2)$ and $(L(f_C),\a_3)$.
Let also $S_{f_A}$, $S_{f_A^{-1}}$, etc., be the kernels giving the corresponding autoequivalences
and their inverses. By Lemma \ref{auteq-conv-lem}, we have 
$$S_L\simeq S_{f_A}\circ_A S_{(f_A\times f_B)(L)}\circ_B S_{f_B^{-1}}[a] \text{ and}$$
$$S_M\simeq S_{f_B}\circ_B S_{(f_B\times f_C)(M)}\circ_C S_{f_C^{-1}}[b]$$
for some $a,b\in\Z$, where we abbreviate $S_{L,\a}$ to $S_L$, etc.. Therefore,
$$S_L\circ_B S_M\simeq S_{f_A}\circ_A S_{(f_A\times f_B)(L)}\circ_B S_{(f_B\times f_C)(M)}\circ_C S_{f_C^{-1}}[a+b].$$
Now Proposition \ref{convolution-prop}(ii) allows to compute the convolution of kernels in the middle, so
we obtain 
$$S_L\circ_B S_M\simeq S_{f_A}\circ_A S_{(f_A\times f_C)(M\circ L)}\circ_C S_{f_C^{-1}}[c]$$
for some $c\in\Z$. By Lemma \ref{auteq-conv-lem}, the right-hand side is isomorphic to
$S_{M\circ L}[i]$ for some $i\in\Z$.
\ed

\begin{cor}\label{Hom-cor}
Let $(L,\a)$ and $(M,\b)$ be generalized Lagrangian pairs for the standard ess-abelian
variety $X_A$. Assume that $L$ and $M$ are transversal, and let $G$ be the corresponding
Heisenberg extension of $L\times_{X_A} M$ (see Proposition \ref{lag-int-prop}), with the
underlying $\G_m$-torsor given by the tensor product of the pull-backs of $\a^{-1}$ and $\b$.
Then $\Hom^*_{D^b(A)}(S_{L,\a},S_{M,\b})$ is concentrated in one degree and is an 
irreducible representation of $G$ of weight one.
\end{cor}

\Pf . We have
$$\Hom^*(S_{L,\a},S_{M,\b})\simeq H^*(A,S_{L,\a}^\vee\ot^{\mathbb L} S_{M,\b}).
$$
Now we interpret $(L,\a)$ and $(M,\b)$ as correspondences from $0$ to $X_A$ and
use Proposition \ref{dual-prop} that says that $S_{L,\a}^\vee$ is isomorphic to a shift
of the LI-kernel corresponding to $(L,\a^{-1})$, viewed as a correspondence from $X_A$ to $0$.
Now we can compute $H^*(A,S_{L,\a}^\vee\ot S_{M,\b})$ by applying Theorem \ref{main-conv-thm}.
Note that the composition of correspondences in our case is $L\times_{X_A} M$, and
$\b\circ\a^{-1}$ is exactly the underlying $\G_m$-torsor of $G$. Thus,
$(M\circ\ov{L},\b\circ\a^{-1})$-invariance structure on a vector space
can be viewed as a weight-$1$ representation of the Heisenberg extension
$G$ (see Example \ref{Lag-ex}.3 and Sec.\ \ref{Heis-gr-sec}).
\ed

\begin{rem} Corollary \ref{Hom-cor} generalizes the well known result that for a nondegenerate
line bundle $L$ on $A$ the cohomology $H^*(A,L)$ is concentrated in one degree and is
an irreducible weight-$1$ representation of the Mumford's theta group attached
to $L$ (see \cite[Sec.\ 16 and 23]{Mum}). Note that actions of some natural groups on the cohomology
of vector bundles on $A$ were also considered by Umemura~\cite{Um}.
\end{rem}

For $g=\left(\begin{matrix} a & b \\ c & d\end{matrix}\right)\in U(X_A,\Q)$ let us denote
$b=b(g)\in\Hom(\hat{A},A)_{\Q}$. 
Let us consider the subset $U^0\sub U(X_A,\Q)$ consisting of $g$ such that
$b(g)$ is invertible. Note that the group $U(X_A,\Q)$ is completely determined by the algebra
$R=\End(A)\ot\Q$ and the Rosati involution $\iota$ on it (with respect to some polarization).
Namely, it consists of automorphisms of the free rank-$2$ module over $R$ preserving
the standard skew-Hermitian form $((x_1,y_1),(x_2,y_2))=\iota(x_1)y_2-\iota(x_2)y_1$.
It follows that $U(X_A,\Q)$ can be identified with the group of $\Q$-points in a connected
algebraic group $U_{X_A}$ (see \cite[Sec.\ 9]{P-thesis} and \cite[Sec.\ 4]{P-maslov}
for a more detailed study of the group $U_{X_A}$, which is denoted there by $\SL_{2,A,\Q}$). 
Since $U(X_A,\Q)$ is Zariski dense 
in $U_{X_A}$ (see \cite[18.3]{Borel}) and since the invertibility of $b(g)$ is a Zariski open
condition, we deduce that
the subset $U^0\sub U(X_A,\Q)$ is {\it big} in the following sense: 
for any triple of elements $g_1,g_2,g_3\in U(X_A,\Q)$
the intersection $(U^0)^{-1}\cap U^0g_1\cap U^0g_2\cap U^0g_3$ is non-empty (this notion
goes back to \cite[IV.\ 42]{Weil} while the term is due to D.~Kazhdan).
The importance of this condition is due to the fact
that a $2$-cocycle of $U(X_A,\Q)$ is uniquely determined by its restriction to $U^0\times U^0$
(see \cite[Lem.\ 4.2]{P-maslov}).

Recall that we denote by $\Lag(X_A)$ the monoid of g-Lag-correspondences $(L,\a)$
from $X_A$ to $X_A$ such that both projections $L\to X_A$ are surjective.
We have a surjective homomorphism $\pi:\Lag(X_A)\to U(X_A,\Q)$ (see Lemma \ref{Lag-Sp-lem}). 
Note that a g-Lag-correspondence $L$ is nondegenerate if and only if $\pi(L)\in U^0$ (see Example 
\ref{nondeg-ex}).
Proposition \ref{convolution-prop} and Theorem \ref{main-conv-thm}
lead to the following computation of the convolution of gLI-kernels
(and hence the composition of gLI-functors)
associated with g-Lag-correspondences from $\Lag(X_A)$.

\begin{thm}\label{Lag-conv-thm} 
(i) For any pair of g-Lag-correspondences $(L,\a),(M,\b)\in\Lag(X_A)$
we have
$$S_{L,\a}\circ_A S_{M,\b}\simeq S_{M\circ L,\b\circ\a}[\la(\pi(M),\pi(L))]$$
for some integer $\la(\pi(L),\pi(M))$. Hence, we also have an isomorphism of functors
$$\Phi_{M,\b}\circ\Phi_{L,\a}\simeq \Phi_{M\circ L,\b\circ\a}[\la(L,M)].$$

\noindent
(ii) The map $\la(g_1,g_2)$ is a $2$-cocycle of $U(X_A,\Q)$ with values in $\Z$.
We have
$$\la(g_1,g_2)=-i(b(g_1)^{-1}b(g_1g_2)b(g_2)^{-1})$$
whenever $b(g_1)$, $b(g_2)$ and $b(g_1g_2)$ are invertible.
Here $i(\cdot)$ denotes the index of a symmetric isogeny (see Section \ref{index-sec}).
\end{thm}

\Pf . (i) Theorem \ref{main-conv-thm} implies that
$$S_{L,\a}\circ_A S_{M,\b}\simeq S_{M\circ L,\b\circ\a}[\la(L,M)]$$
for some integer $\la(L,M)$. 
We have to prove that $\la(L,M)$ depends only on $\pi(L)$ and $\pi(M)$.
Let $g=\pi(L)$ and let us equip $L(g)\sub X_A\times X_A$ with a line bundle $\a_g$
so that $(L(g),\a_g)$ is a Lagrangian correspondence. By Proposition \ref{lag-inv-prop}(ii),
$S_{L,\a}$ is obtained by successive extensions from objects of the form
$T_{x}(S_{L(g),\a_g})$
where $x\in X_{A\times A}$. Using Lemmas  \ref{T-A-B-lem} and \ref{K-T-lem},
we can rewrite such objects as $\bbT_{x'}\circ_A S_{L(g),\a_g}$
with $x'\in X_A$. This immediately implies that $S_{L,\a}\circ_A S_{M,\b}$
is obtained by successive extensions from objects of the form
$\bbT_{x'}\circ_A S_{L(g),\a_g}\circ_A S_{M,\b}$, hence
$\la(L,M)=\la(L(g),M)$. A similar argument shows that $\la(L,M)$ depends only on $\pi(M)$.

\noindent
(ii) The fact that $\la(\cdot,\cdot)$ is a $2$-cocycle follows from the definition.
Let us show how to rewrite the formula \eqref{la-L-M-eq} in the required form
for $L=L(g_2)$ and $M=L(g_1)$, where $g_1,g_2\in U^0$ are such that $g_1g_2\in U^0$. Let
$$g_i=\left(\begin{matrix} \a_i & \b_i \\ \ga_i & \de_i\end{matrix}\right), \ i=1,2,$$
where $\b_1$ and $\b_2$ are invertible.
For $i=1,2$ we have isomorphisms in $\Ab_\Q$
$$\phi_i: A\times\hat{A}\to L(g_i): (a,\xi)\mapsto (a, \xi, \a_i a+\b_i\xi, \ga_i a+ \de_i\xi),
$$
Recall that $-\la(M,L)$ is the index of the symmetric isogeny
$\hat{f_0}g_0:F_0\to\hat{F_0}$, where
$F=\ker(L(g_2)\times_A L(g_1)\rTo{p_{17}} A\times A)$
and $f_0$ and $g_0$ were defined in Proposition \ref{convolution-prop}(ii).
Here we view $L(g_2)\times_A L(g_1)$ as a subvariety in 
$X_A\times X_A\times X_A\times X_A=A\times\hat{A}\times\ldots $ (the last product has
$8$ factors) and denote by $p_{17}$ the corresponding projection to $A\times A$.
It is easy to check that there is an isomorphism in $\Ab_\Q$
$$\phi:\hat{A}\to F_0:\xi\mapsto (0, \xi, \b_2\xi, \de_2\xi), 
(\b_2\xi, -\ga\xi, 0, (\ga_1\b_2-\de_1\ga)\xi),$$
where $\ga=\b_1^{-1}\a_1\b_2$.
Furthermore, we have
$$f_0(\phi(\xi))=\b_2\xi \text{ and }\ g_0(\phi(\xi))=(\de_2+\ga)\xi,$$
Hence,
$$\hat{\phi}\hat{f}g\phi=\hat{\b}_2(\de_2+\ga)=\hat{\b}_2(\de_2+\b_1^{-1}\a_1\b_2)$$
and we obtain
$$\la(g_1,g_2)=i(\hat{\b}_2(\de_2+\b_1^{-1}\a_1\b_2))=i((\de_2+\b_1^{-1}\a_1\b_2)\b_2^{-1})=
i(\de_2\b_2^{-1}+\b_1^{-1}\a_1)$$
which is equivalent to the desired formula.
\ed

\begin{rem}
In the case when $g_1,g_2\in U(X_A,\Q)\cap\End(X_A)$, our formula for $\la(g_1,g_2)$
agrees with that of Orlov in \cite[Sec.\ 4]{O}, due to the standard formula for the index
of a line bundle (see Sec.\ \ref{index-sec}).
\end{rem}

\subsection{Central extensions related to LI-endofunctors}\label{central-ext-sec}

From now on we assume that $\cha(k)=0$.

\begin{defi}
Let us denote by $\KER(A,A)$ the set of isomorphism classes of objects of $D^b(A\times A)$.
We equip it with a semiring structure by taking the direct sum $\oplus$ as addition
and the convolution $\circ_A$ (see Sec.\ \ref{conv-sec}) as multiplication.
\end{defi}

Note that we have a homomorphism of semirings
$$\KER(A,A)^{op}\to\Fun(D^b(A),D^b(A)):K\mapsto \Phi_K,$$
where $\Phi_K$ is the Fourier-Mukai functor given by the kernel $K$ (see Sec.\ \ref{conv-sec})
and $\Fun(\cdot,\cdot)$ is the set of isomorphism classes of exact functors .


\begin{defi}
(i) Let us say that a cohomologically pure kernel $K\in\KER(A,A)$ is an
{\it LI-kernel} (resp., {\it weak LI-kernel}) if
there exists a Lagrangian correspondence $(L,\a)$ in $\Lag(X_A)$ such that
$K$ is $(L,\a)$-invariant (resp., $K=\bigoplus_{i=1}^N K_i$, where $K_i$
are $(L,\a_i)$-invariant for some Lagrangian correspondences $(L,\a_i)$ with common $L$).

\noindent
(ii) Let $\KER^{LI}(A,A)$ (resp., $\KER^{wLI}(A,A)\sub\KER(A,A)$; resp.,
$\KER^{LI}_{\oplus}(A,A)\sub\KER(A,A)$)
denote the subset consisting of LI-kernels (resp., weak LI-kernels; 
resp., finite direct sums of LI-kernels).
By Theorem \ref{Lag-conv-thm}(i) and Corollary \ref{gen-Lag-cor}
$\KER^{LI}_{\oplus}(A,A)$ is a subsemiring in $\KER(A,A)$, while
$\KER^{wLI}(A,A)$ is a submonoid in the mutliplicative monoid $(\KER(A,A),\circ_A)$.

\noindent
(iii)
Consider the following equivalence relation on $\KER^{LI}(A,A)$ (resp., $\KER^{wLI}(A,A)$):
$K\sim_{\bH} K'$ if there exists a point $x\in X_{A\times A}$ such that $K'\simeq T_x(K)$
(resp., there exists a direct sum decomposition $K=\bigoplus_{i=1}^N K_i$
and a collection of points $x_i\in X_{A\times A}$, $i=1,\ldots,N$, such that 
$K'\simeq \bigoplus_{i=1}^N T_{x_i}(K_i)$.)
By Lemmas \ref{T-A-B-lem} and \ref{K-T-lem}, this is equivalent to the existence
of $x'\in X_A$ and $x''\in X_A$ such that 
$$K'\simeq \bbT^A_{x'}\circ_A K\simeq K\circ_A \bbT^A_{x''}.$$
This easily implies that the set of equivalence classes 
$$\ov{\KER}^{LI}(A,A):=\KER^{LI}(A,A)/\sim_{\bH}\ \ \simeq\ \KER^{wLI}(A,A)/\sim_{\bH}$$
inherits the monoid structure (with respect to the convolution $\circ_A$).
\end{defi}

We have a central submonoid 
$Z\sub \ov{\KER}^{LI}(A,A)$ consisting of kernels of the form
$\OO_{\De_A}^{\oplus m}[n]$, where $\De_A\sub A\times A$
is the diagonal. The convolution with an element in $Z$ is a very simple operation:
$$K\circ_A \OO_{\De}^{\oplus m}[n]\simeq K^{\oplus m}[n].$$
We will use the notation $m\cdot K:= K^{\oplus m}$.
Note that $Z$ is isomorphic to $\N^*\times\Z$, so the localization of $Z$ with respect
to the multiplicative set $\N^*\sub Z$ gives the group $Z_{\Q^*}=\Q^*\times\Z$.


Let us denote by $\ov{\KER}^{LI}(A,A)_{\Q^*}$ the localization of $\ov{\KER}^{LI}(A,A)$
with respect to the central multiplicative set $\N^*\sub Z$.

\begin{lem}
The monoid $\ov{\KER}^{LI}(A,A)_{\Q^*}$ is a group.
\end{lem}

\Pf . Let $(L,\a)$ be a Lagrangian correspondence from $X_A$ to $X_A$, such that
both projections $L\to X_A$ are surjective. Consider the transposed Lagrangian
correspondence $(\si L,\a^{-1})$ and the composed g-Lag-correspondence
$(Z,\b)=(L\circ \si L, \a\circ\a^{-1})$. Let $j:Z\to X_A\times X_A$ be the natural map, and let
$Z_0$ be the connected component of zero in $Z$.
Since we have the diagonal map $L\to Z$, we obtain that $j(Z_0)$ is equal to the
diagonal $\De_X\sub X_A\times X_A$. Therefore, by Theorem \ref{Lag-conv-thm} and
Proposition \ref{lag-inv-prop}(ii), we obtain
$$S_{L,\a}\circ_A S_{\si L,\a^{-1}}=N\cdot S_{\De_X,\OO}[i]=N\cdot \OO_{\De_A}[i]$$
in $\ov{\KER}^{LI}(A,A)$
with $N\in \N^*$ and $i\in\Z$. This shows that $S_{L,\a}$ has a right inverse in
$\ov{\KER}^{LI}(A,A)_{\Q^*}$. Exchanging the roles of $L$ and $\si L$ we also obtain
the existence of a left inverse. Since every element
of $\KER^{LI}(A,A)$ is of the form $n\cdot S_{L,\a}$ with $n\in\N^*$, we see that
it is invertible in $\ov{\KER}^{LI}(A,A)_{\Q^*}$.
\ed

The first two parts of the following theorem summarize most of the picture discussed above.

\begin{thm} Assume that $\cha(k)=0$.

\noindent
(i) The semiring $\KER^{LI}_{\oplus}(A,A)$ has as an additive basis
the kernels of the form $S_{L(g),\a}[n]$, where 
$L(g)$ is the Lagrangian correspondence from $X_A$ to $X_A$ associated with $g\in U(X_A,\Q)$
(see Example \ref{Lag-corr-ex}), 
$\a$ is a line bundle on $L(g)$ such that $(L(g),\a)$ is a Lagrangian pair, and $n$ is
an integer. The multiplication in $\KER^{LI}_{\oplus}(A,A)$ has form
\begin{equation}\label{semiring-eq}
S_{L(g_2),\a_2}[n_2]\circ_A S_{L(g_1),\a_1}[n_1]= 
\bigoplus_{i=1}^N S_{L(g_1g_2),\b_i}[n_1+n_2+\la(g_1,g_2)]
\end{equation}
with $\la(\cdot,\cdot)$ as in Theorem \ref{Lag-conv-thm}
and some line bundles $\b_1,\ldots,\b_N$ on $L(g_1g_2)$, such that $(L(g_1g_2),\b_i)$
are Lagrangian pairs. 

\noindent
(ii) For each $g\in U(X_A,\Q)$ let us choose a line bundle $\a_g$ in such a way that $(L(g),\a_g)$ is a Lagrangian pair. Then the class $S(g)$ of the kernel 
$S_{L(g),\a_g}\in\KER^{LI}(A,A)$ in $\ov{\KER}^{LI}(A,A)$ depends only on $g$.
For $g_1,g_2\in U(X_A,\Q)$ 
one has the following equality in $\ov{\KER}^{LI}(A,A)$:
\begin{equation}\label{S-comp-eq}
S(g_2)\circ_A S(g_1)=N(g_1,g_2)\cdot S(g_1g_2)[\la(g_1,g_2)],
\end{equation}
where $N(g_1,g_2)\in\N^*$.
The map $g\mapsto S(g)$
extends to an isomorphism of groups
\begin{equation}\label{main-homomorphism-eq}
\widehat{U}(X_A,\Q)\to \ov{\KER}^{LI}(A,A)_{\Q^*}^{op},
\end{equation}
where $\widehat{U}(X_A,\Q)$ is the central extension of $U(X_A,\Q)$ by $Z_{\Q^*}=\Q^*\times\Z$
associated with the $2$-cocycle $(N(g_1,g_2),\la(g_1,g_2))$.

\noindent (iii)
For $g_1,g_2\in U(X_A,\Q)$ 
such that $b(g_1)$, $b(g_2)$ and $b(g_1g_2)$ are invertible, one has
\begin{equation}\label{N-g1-g2-eq}
N(g_1,g_2)=\frac{q(L(g_1))^{1/2}q(L(g_2))^{1/2}}{q(L(g_1g_2))^{1/2}}\in \N^*.
\end{equation}
with $q(L)$ given by \eqref{q-L-eq}.
\end{thm}

\Pf . The assertion about the additive basis in $\KER^{LI}_{\oplus}(A,A)$ follows
from Theorem \ref{lag-inv-thm} 
and Corollary \ref{uniqueness-cor}. The equation \eqref{semiring-eq} follows from
Theorem \ref{Lag-conv-thm} together with Corollary \ref{gen-Lag-cor}.

The fact that the class of $S_{L(g),\a_g}$
in  $\ov{\KER}^{LI}(A,A)$ does not depend on the choice of $\a_g$ follows immediately 
from Lemma \ref{beta-change-lem}. Furthermore, together with Corollary
\ref{uniqueness-cor} this implies that the map $(g,m,n)\to S(g)^{\oplus m}[n]$ is a bijection
between $U(X_A,\Q)\times\N^*\times\Z$ and $\ov{\KER}^{LI}(A,A)$.

Applying Proposition \ref{lag-inv-prop}(ii) to $Z=L(g_1)\circ L(g_2)$ we get
$$S_{L(g_1)\circ L(g_2),\a_{g_1}\circ\a_{g_2}}\simeq
\left(\bigoplus_{x\in \Pi'} T_x(S_{L(g_1g_2),\a_{g_1g_2}})\right)^{\oplus N_Z}.$$
Recall that here
$$|\Pi'|=|\pi_0(j(Z))|=d:=\deg(Z_0\to j(Z_0))$$ 
(see \eqref{C-Pi-eq}) and
$N_Z=\frac{|\pi_0(Z)|^{1/2}}{d^{1/2}}$.
Hence, by Theorem \ref{Lag-conv-thm}, in $\ov{\KER}^{LI}(A,A)$ we have
$$S_{L(g_1)\circ L(g_2),\a_{g_1}\circ\a_{g_2}}=|\pi_0(Z)|^{1/2}\cdot d^{1/2}
\cdot S(g_1g_2)$$
which implies \eqref{S-comp-eq} with
$$N(g_1,g_2)=|\pi_0(Z)|^{1/2}\cdot d^{1/2}.$$
This gives a homomorphism \eqref{main-homomorphism-eq} from the central extension
$\widehat{U}(X_A,\Q)$, which is easily seen to be an isomorphism from the above identification
of the set $\ov{\KER}^{LI}(A,A)$.

It remains to prove the formula \eqref{N-g1-g2-eq} in the case when
$b(g_1)$, $b(g_2)$ and $b(g_1g_2)$ are invertible. 
Note that by Lemma \ref{q-lem}, we have
$$q(Z)=q(L(g_1))\cdot q(L(g_2)).$$
Thus, \eqref{N-g1-g2-eq} can be rewritten as
$$\pi_0(Z)\cdot d=\frac{q(Z)}{q(j(Z_0))}=\frac{\deg(Z\to X_A)}{\deg(j(Z_0)\to X_A)}.$$
But this can be checked by comparing the degrees in the commutative diagram of isogenies
\begin{diagram}
Z_0 &\rTo{}& Z\\
\dTo{j}&&\dTo{}\\
j(Z_0) &\rTo{}& X
\end{diagram}
as in the end of proof of Proposition \ref{lag-inv-prop}(ii).
\ed

\begin{cor}
Let $\ov{\KER}^{LI}(A,A)_{\R^*}$ be the push-out of the central
extension sequence
$$1\to \Q^*\to \ov{\KER}^{LI}(A,A)_{\Q^*}\to \ov{\KER}^{LI}(A,A)_{\Q^*}/\Q^*\to 1$$
with respect to the embedding $\Q^*\to\R^*$, and let
$\wt{U}(X,\Q)$ be the central extension of $U(X,\Q)$ by $\Z$ associated with the 2-cocycle
$\la(\cdot,\cdot)$. Then the map $g\mapsto S(g)\cdot q(L(g))^{-1/2}$ extends to
a homomorphism
$$\wt{U}(X,\Q)\to \ov{\KER}^{LI}(A,A)_{\R^*}.$$
\end{cor}

The $2$-cocycle $N(\cdot,\cdot)$ given by \eqref{N-g1-g2-eq} gives in general a nontrivial
cohomology class in $H^2(U(X_A,\Q),\Q^*)$ as the following result shows.

\begin{prop}
Let $E$ be an elliptic curve with complex multiplication by $\sqrt{D}$, where $D<-1$ and
$|D|$ is square-free.
Then the cohomology class in $H^2(U(X_E,\Q),\Q^*)$ given by $N(\cdot,\cdot)$ is nontrivial.
\end{prop}

\Pf . Let $\sqrt{\Q^*}\sub\R^*$ be the subgroup of $x\in\R^*$ such that $x^2\in\Q^*$.
We have a homomorphism
$$\ov{q}^{1/2}:U(X_E,\Q)\to\sqrt{\Q^*}/\Q^*: g\mapsto q(L(g))^{1/2} \mod \Q^*$$
that we can view as an element of $H^1(U(X_E,\Q),\sqrt{\Q^*}/\Q^*)$.
Now the class of $N(\cdot,\cdot)$ is equal to $\de(q^{1/2})$ where
$$\de:H^1(U(X_E,\Q),\sqrt{\Q^*}/\Q^*)\to H^2(U(X_E,\Q),\Q^*)$$
is the coboundary map coming from the exact sequence
$$1\to \Q^*\to \sqrt{\Q^*}\to \sqrt{\Q^*}/\Q^*\to 1.$$
Thus, it is enough to prove that $\ov{q}^{1/2}$ cannot be lifted to a homomorphism
$U(X_E,\Q)\to\sqrt{\Q^*}$. We are going to prove that it cannot even be lifted to a homomorphism
$U(X_E,\Q)\to\sqrt{\Q^*}/\pm 1$. Indeed, using the isomorphism 
$\sqrt{\Q^*}/\pm 1\to \Q^*:x\mapsto x^2$ we see that this is equivalent to asserting that
the homomorphism
$$\ov{q}:Sp(X_E,\Q)\to\Q^*/(\Q^*)^2: g\mapsto q(L(g)) \mod (\Q^*)^2$$
cannot be lifted to a homomorphism to $\Q^*$. 

Note that the group $U(X_E,\Q)$ in our case is isomorphic to the unitary group
$U=U_2(K,f)$, where $K=\Q(\sqrt{D})$ and $f$ is a skew-Hermitian form of index $1$
(we use the terminology of \cite{Dieud-cl-gr}). Let us consider the natural embedding 
$$\iota:K^*/\Q^*\to K^*: a\mapsto \frac{\ov{a}}{a}.$$
It is known that the image of the determinant map $\det:U\to K^*$ is contained
in $\iota(K^*/\Q^*)$ (see \cite[Thm.\ 3]{Dieud}). Hence, we have a unique homomorphism
$\phi:U\to K^*/\Q^*$ such that $\det(g)=\iota(\phi(g))$.
We claim that
$$\ov{q}(g)=\ov{\Nm}(\phi(g)),$$
where $\ov{\Nm}:K^*/\Q^*\to \Q^*/(\Q^*)^2$ is induced by the norm homomorphism
$\Nm:K^*\to \Q^*$. Indeed, it is known that in our case
the special unitary group $SU\sub U$ coincides with the normal subgroup $T\sub U$ generated
by unitary transvections (see \cite[\S 5]{Dieud-cl-gr}), and the latter subgroup coincides with 
$\SL_2(\Q)\sub U$. It is easy to see that for $g\in\SL_2(\Q)$ the degree of the
projection $p_{12}:L(g)\to E^2$ is a square.
Thus, $\ov{q}$ is trivial on $SU$ and we should have
$$\ov{q}(g)=\chi(\phi(g))$$
for some homomorphism $\chi:K^*/\Q^*\to \Q^*/(\Q^*)^2$. It remains to check that our
statement holds for diagonal matrices
$$g_a=\left(\begin{matrix} a^{-1} & 0 \\ 0 & \ov{a} \end{matrix}\right)$$
where $a\in K^*$ (since $\phi(g_a)=a$). In other words, we claim that in this case
$$\ov{q}(L(g_a))\equiv\Nm(a)\mod (\Q^*)^2.$$
It is enough to check this in the case when $a\in\End(E)\sub K$.
Then the Lagrangian correspondence $L(g_a)$ can be described as
$$L(g_a)=\{(ax,y,x,\ov{a}y)\ |\ x,y\in E\}\sub E^2\times E^2,$$
so $q(L(g_a))=\deg(a)=\Nm(a)$.

Finally, suppose we have a homomorphism $\rho:U\to \Q^*$ such that
$\ov{\Nm}(\phi(g))=\rho(g)\mod (\Q^*)^2$. 
Let us specialize this equality to $g=g_{\sqrt{D}}$.
Then we should have 
$$\rho(g_{\sqrt{D}})\equiv -D\mod (\Q^*)^2.$$
But $\rho$ factors through $U/[U,U]$ and $g_{\sqrt{D}}$ projects to
an element of finite order in $U/[U,U]$. Therefore, $\rho(g_{\sqrt{D}})$ should
be an element of finite order in $\Q^*$, i.e., $\rho(g_{\sqrt{D}})=\pm 1$ which
is a contradiction.
\ed


\end{document}